\documentclass[MSNbibl,number,citesort,dvips]{arxstspdf}
\usepackage{subenv}
\usepackage{graphicx}
\usepackage{flushend}
\usepackage{stfloats}

\volume{27}
\issue{4}
\pubyear{2012}
\firstpage{538}
\lastpage{557}
\doi{10.1214/12-STS400} %kopijuoti is PTS
\makeatletter

\newtheorem{theos}{Theorem}
\newtheorem{props}{Proposition}
\newtheorem{lems}{Lemma}
\newtheorem{cors}{Corollary}

\newproclaim{exas}{Example}
\newproclaim{defns}{Definition}
\newproclaim{Remarks}{Remarks}

\newcommand{\norm}{|\!|\!|}
\newcommand{\llangle}{\langle\!\langle}
\newcommand{\rrangle}{\rangle\!\rangle}

\newcommand{\rscslop}{\tau_\Loss}
\newcommand{\taupar}{{\eta}}

\newcommand{\Exs}{{{\mathbb{E}}}}

\newcommand{\trace}{\operatorname{trace}}
\newcommand{\sign}{\operatorname{sign}}

\newcommand{\betapar}{{\theta}}
\newcommand{\numobs}{{n}}
\newcommand{\pdim}{{p}}
\newcommand{\kdim}{{s}}
\newcommand{\regpar}{{\lambda_\numobs}}
\newcommand{\Sset}{{S}}

\newcommand{\Ball}{{\mathbb{B}}}
\newcommand{\Loss}{{\mathcal{L}}}

\newcommand{\Regplain}{{\mathcal{R}}}

\newcommand{\ConvSet}{{\mathbb{C}}}

\newcommand{\Mcon}{{\kappa}}
\newcommand{\NewMcon}{{\Mcon_{\Loss}}}

\newcommand{\qpar}{{q}}
\newcommand{\delpar}{{\Delta}}
\newcommand{\mprob}{{\mathbb{P}}}
\newcommand{\Xmat}{{X}}

\newcommand{\radq}{{R_\qpar}}
\newcommand{\param}{\theta}

\newcommand{\real}{\mathbb{R}}
\newcommand{\plaincon}{{c}}

\newcommand{\TmpVarA}{{\theta}}
\newcommand{\TmpVarB}{{\gamma}}

\newcommand{\ConeSet}{{\mathbb{C}}}

\newcommand{\rowsp}{\operatorname{row}}
\newcommand{\colsp}{\operatorname{col}}

\newcommand{\ModelSet}{{\mathcal{M}}}

\newcommand{\pdima}{{{\pdim_1}}}
\newcommand{\pdimb}{{{\pdim_2}}}
\newcommand{\Apperr}{{\mathcal{E}_{\mathrm{app}}}}
\newcommand{\Esterr}{{\mathcal{E}_{\mathrm{err}}}}

\newcommand{\Stau}{{\Sset_\taupar}}

\newcommand{\delLike}{{\delta\Loss}}

\newcommand{\numgroup}{{{N_{\GROUP}}}}

\newcommand{\gsize}{{m}}

\newcommand{\Gmax}{{m}}

\newcommand{\gqpar}{{\alpha}}
\newcommand{\gqvec}{{\gqpar}}

\newcommand{\GROUP}{{\mathcal{G}}}
\newcommand{\DualNoise}{{\rho_{\GROUP}}}

\newcommand{\basgauss}{{\varepsilon}}

\newcommand{\GroupSset}{{\Sset_{\GROUP}}}
\newcommand{\groupcard}{{\kdim_\GROUP}}

\newcommand{\delparhatlam}{{{\widehat{\delpar}}_{\regpar}}}

\newcommand{\SUPERGVEC}{{\vec{\gqpar}}}
\newcommand{\gqparbar}{{\gqpar}}

\newcommand{\Risk}{{\overline{\Loss}}}

\newcommand{\curve}{{\kappa}}

\newcommand{\GARTWO}{{\kappa_2}}

\makeatother

\begin{document}
\begin{frontmatter}

\title{A Unified Framework for High-Dimensional Analysis of
$M$-Estimators with Decomposable Regularizers}%\thanksref{T1}
% kai straipsnis turi susijusiu diskusiju ir rejoinder'iu
%rejoinder at \relateddoi{r}{10.1214/00-STSXXXX}.}
\runtitle{High-Dimensional Analysis of Regularized $M$-Estimators}

\begin{aug}
\author[a]{\fnms{Sahand N.} \snm{Negahban}\ead[label=e1]{sahandn@mit.edu}},
\author[b]{\fnms{Pradeep} \snm{Ravikumar}\ead[label=e2]{pradeepr@cs.utexas.edu}},
\author[c]{\fnms{Martin J.} \snm{Wainwright}\corref{}\ead[label=e3]{wainwrig@stat.berkeley.edu}}
\and
\author[d]{\fnms{Bin} \snm{Yu}\ead[label=e4]{binyu@stat.berkeley.edu}}
\runauthor{Negahban, Ravikumar, Wainwright and Yu}

\affiliation{Massachusetts Institute of Technology, University of
Texas, University of California and University of California}

\address[a]{Sahand Negahban is Postdoctoral Researcher,
EECS Department, Massachusetts Institute of Technology,
\mbox{77 Massachusetts} Avenue, Cambridge, Massachusetts 02139, USA
\printead{e1}.}
\address[b]{Pradeep Ravikumar is Assistant Professor, Department of Computer
Science, University of Texas, Austin, Texas 78712, USA \printead{e2}.}
\address[c]{Martin J. Wainwright is Professor, Departments of Statistics and EECS,
University of California, Berkeley, California 94720, USA
\printead{e3}.}
\address[d]{Bin Yu is Professor, Departments of Statistics and EECS, University of
California, Berkeley, California 94720, USA \printead{e4}.}

\end{aug}

% ABSTRACT
%
\begin{abstract}
High-dimensional statistical inference deals with models in which the
the number of parameters $p$ is comparable to or larger than the
sample size $n$. Since it is usually impossible to obtain consistent
procedures unless $p/n \rightarrow0$, a line of recent work has
studied models with various types of low-dimensional structure,
including sparse vectors, sparse and structured matrices, low-rank
matrices and combinations thereof. In such settings, a general
approach to estimation is to solve a regularized optimization problem,
which combines a loss function measuring how well the model fits the
data with some regularization function that encourages the assumed
structure. This paper provides a unified framework for establishing
consistency and convergence rates for such regularized $M$-estimators
under high-dimensional scaling. We state one main theorem and show
how it can be used to re-derive some existing results, and also to
obtain a number of new results on consistency and convergence rates,
in both $\ell_2$-error and related norms. Our analysis also identifies
two key properties of loss and regularization functions, referred to
as restricted strong convexity and decomposability, that ensure
corresponding regularized $M$-estimators have fast convergence rates
and which are optimal in many well-studied cases.
\end{abstract}

% KEYWORDS
%

\end{frontmatter}

%s1 #&#
\section{Introduction}

High-dimensional statistics is concerned with models in which the
ambient dimension of the problem $\pdim$ may be of the same order
as---or substantially larger than---the sample size $\numobs$. On the
one hand, its roots are quite old, dating back to work on random matrix
theory and high-dimensional testing problems (e.g., [\cite*{Gir95},
\cite*{Mehta}, \cite*{Pas72,Wig55}]). On the other hand, the past
decade has witnessed a tremendous surge of research activity. Rapid
development of data collection technology is a major driving force: it
allows for more observations to be collected (larger $\numobs$) and
also for more variables to be measured (larger $\pdim$). Examples are
ubiquitous throughout science: astronomical projects such as the Large
Synoptic Survey Telescope (available at
\href{http://www.lsst.org/lsst/}{www.lsst.org}) produce terabytes of
data in a single evening; each sample is a high-resolution image, with
several hundred megapixels, so that $\pdim\gg10^8$. Financial data is
also of a high-dimensional nature, with hundreds or thousands of
financial instruments being measured and tracked over time, often at
very fine time intervals for use in high frequency trading. Advances in
biotechnology now allow for measurements of thousands of genes or
proteins, and lead to numerous statistical challenges (e.g., see the
paper \cite{BicBroHuaLi09} and references therein). Various types of
imaging technology, among them magnetic resonance imaging in medicine
\cite{LusDonSanPau08} and hyper-spectral imaging in
ecology~\cite{Lan02}, also lead to high-dimensional data sets.

In the regime $\pdim\gg\numobs$, it is well known that consistent
estimators cannot be obtained unless additional constraints are
imposed on the model. Accordingly, there are now several lines of
work within high-dimensional statistics, all of which are based on
imposing some type of low-dimensional constraint on the model space
and then studying the behavior of different estimators. Examples
include linear regression with sparsity constraints, estimation of
structured covariance or inverse covariance matrices, graphical model
selection, sparse principal component analysis, low-rank matrix
estimation, matrix decomposition problems and estimation of sparse
additive nonparametric models. The classical technique of
regularization has proven fruitful in all of these contexts. Many
well-known estimators are based on solving a convex optimization
problem formed by the sum of a loss function with a weighted
regularizer; we refer to any such method as a \textit{regularized
\mbox{$M$-estimator}}. For instance, in application to linear models, the
Lasso or basis pursuit approach \cite{Tibshirani96,Chen98} is based
on a combination of the least squares loss with
$\ell_1$-regularization, and so involves solving a quadratic program.
Similar approaches have been applied to generalized linear models,
resulting in more general (nonquadratic) convex programs with
$\ell_1$-constraints. Several types of regularization have been used
for estimating matrices, including standard $\ell_1$-regularization, a
wide range of sparse group-structured regularizers, as well as
regularization based on the nuclear norm (sum of singular values).

%pa1.subsection.subsubsection.1 #&#
\subsection*{Past Work} Within the framework of
high-dimensional statistics, the goal is to obtain bounds on a given
performance metric that hold with high probability for a finite sample
size, and provide explicit control on the ambient dimension $\pdim$,
as well as other structural parameters such as the sparsity of a
vector, degree of a graph or rank of matrix. Typically, such bounds
show that the ambient dimension and structural parameters can grow as
some function of the sample size $\numobs$, while still having the
statistical error decrease to zero. The choice of performance metric
is application-dependent; some examples include prediction error,
parameter estimation error and model selection error.

By now, there are a large number of theoretical results\vspace*{1pt} in
place for various types of regularized $M$-estimators.\footnote{Given
the extraordinary number of papers that have appeared in recent years,
it must be emphasized that our referencing is necessarily incomplete.}
Sparse linear regression has perhaps been the most active area, and
multiple bodies of work can be differentiated by the error metric under
consideration. They include work on exact recovery for noiseless
observations (e.g., \cite {DonTan05,Donoho06,CandesTao05}), prediction
error consistency (e.g., \cite{GreenRitov04,Bun07,GeerBuhl09,Huang06}),
consistency of the parameter estimates in $\ell_2$ or some other norm
(e.g., \cite
{Bun07,BunWegTsyb07,GeerBuhl09,Huang06,MeiYu09,BiRiTsy08,CandesTao06}),
as well as variable selection consistency (e.g.,
\cite{Meinshausen06,Wainwright06,Zhao06}). The information-theoretic
limits of sparse linear regression are also well understood, and
$\ell_1$-based methods are known to be optimal for $\ell_\qpar$-ball
sparsity \cite{RasWaiYu09} and near-optimal for model selection
\cite{Wainwright09info}. For generalized linear models (GLMs),
estimators based on $\ell_1$-regularized maximum likelihood have also
been studied, including results on risk consistency \cite{Geer08},
consistency in the $\ell_2$ or $\ell_1$-norm \cite{Bach10,Kak10,Mei08}
and model selection consistency \cite{RavWaiLaf08,Bunea08}. Sparsity
has also proven useful in application to different types of matrix
estimation problems, among them banded and sparse covariance matrices
(e.g., \mbox{\cite{BicLev08b,CaiZhou10,ElKaroui}}). Another line of work has
\mbox{studied} the problem of estimating Gaussian Markov random fields or,
equivalently, inverse covariance matrices with sparsity constraints.
Here there are a range of results, including convergence rates in
Frobenius, operator and other matrix norms
\cite{Rot09,RavWaiRasYu2011,LamFan07,ZhoLafWas08},\vadjust{\goodbreak} as well as results
on model selection consistency
\cite{RavWaiRasYu2011,LamFan07,Meinshausen06}. Motivated by
applications in which sparsity arises in a structured manner, other
researchers have proposed different types of block-structured
regularizers (e.g.,
\cite{Tro06,Kim06,Turlach05,ZhaRoc06,YuaLi06,Bach08a,Bar08,Jac09}),
among them the group Lasso based on $\ell_1/\ell_2$-regularization.
High-dimensional consisten\-cy results have been obtained for exact
recovery based on noiseless observations \cite{StohParHas09,Bar08},
convergence rates in the $\ell_2$-norm (e.g.,
\cite{NarRin08,HuaZha09,Lou09,Bar08}) as well as model selection
consistency (e.g., \mbox{\cite{OboWaiJor08,NegWai08,NarRin08}}). Problems of
low-rank matrix estimation also arise in numerous applications.
Techniques based on nuclear norm regularization have been studied for
different statistical models, including compressed sensing
\cite{RecFazPar10,LeeBres09}, matrix completion
\cite{CanRec08,Keshavan09b,Recht09,NegWai10b}, multitask regression
\mbox{\cite{YuaEkiLuMon07,NegWai09,RohTsy10,Bunea10,Bach08b}} and system
identification \mbox{\cite{Fa02,NegWai09,LiuVan09}}. Finally, although the
primary emphasis of this paper is on high-dimensional parametric
models, regularization methods have also proven effective for a class
of high-dimensional nonparametric models that have a sparse additive
decomposition (e.g., \cite{RavLiuLafWas08,Meier09,KolYua08,KolYua10}),
and have been shown to achieve minimax-optimal rates
\cite{RasWaiYu10b}.

%pa1.subsection.subsubsection.2 #&#
\subsection*{Our Contributions} As we have noted previously,
almost all of these estimators can be seen as particular types of
regularized $M$-estimators, with the choice of loss function,
regularizer and statistical assumptions changing according to the
model. This methodological similarity suggests an intriguing
possibility: is there a \textit{common set of theoretical principles}
that underlies analysis of all these estimators? If so, it could be
possible to gain a unified understanding of a large collection of
techniques for high-dimensional estimation and afford some insight
into the literature.

The main contribution of this paper is to provide an affirmative
answer to this question. In particular, we isolate and highlight two
key properties of a regularized $M$-estimator---namely, a
\textit{decomposability property} for the regularizer and a notion of
\textit{restricted strong convexity} that depends on the interaction
between the regularizer and the loss function. For loss functions and
regularizers satisfying these two conditions, we prove a general
result (Theorem~\ref{ThmMain}) about consistency and convergence rates
for the associated estimators. This result provides a family of
bounds indexed by subspaces, and each bound consists of the sum of
approximation error and estimation error. This general result, when
specialized to different statistical models, yields in a direct manner
a large number of corollaries, some of them known and others novel.
In concurrent work, a subset of the current authors has also used
this framework to prove several results on low-rank matrix estimation
using the nuclear norm \cite{NegWai09}, as well as minimax-optimal
rates for noisy matrix completion \cite{NegWai10b} and noisy matrix
decomposition \cite{AgaNegWai11}. Finally, en route to establishing
these corollaries, we also prove some new technical results that are
of independent interest, including guarantees of restricted strong
convexity for group-structured regularization
(Proposition~\ref{PropREGroup}).

The remainder of this paper is organized as follows. We begin in
Section~\ref{SecBackground} by formulating the class of regularized
$M$-estimators that we consider, and then defining the notions of
decomposability and restricted strong convexity.
Section~\ref{SecMain} is devoted to the statement of our main result
(Theorem~\ref{ThmMain}) and discussion of its consequences.
Subsequent sections are devoted to corollaries of this main result for
different statistical models, including sparse linear regression
(Section~\ref{SecSparseReg}) and estimators based on group-structured
regularizers (Section~\ref{SecGroup}). A number of technical results
are presented within the appendices in the supplementary
file \cite{Neg12supp}.
%

%s2 #&#
\section{Problem Formulation and Some Key Properties}
\label{SecBackground}

In this section we begin with a precise formulation of the problem,
and then develop some key properties of the regularizer and loss
function.

%s2.1 #&#
\subsection{A Family of $M$-Estimators}

Let ${Z_1^{n}}:=\{Z_1,\ldots, Z_\numobs\}$ denote
$\numobs$
identically distributed observations with marginal distribution
$\mprob$, and suppose that we are interested in estimating some
parameter $\theta$ of the distribution $\mprob$. Let $\Loss\dvtx
{\real^\pdim}\times\mathcal{Z}^\numobs\rightarrow\real$ be a convex
and differentiable loss function that, for a given set of observations
${Z_1^{\numobs}}$, assigns a cost $\Loss(\betapar; {Z_1^{n}})$ to any
parameter $\param\in\real^\pdim$. Let $\param^*\in\arg\min_{\theta\in
\real^{\pdim}}\Risk(\theta)$ be any minimizer of
the population risk $\Risk(\theta) :=
\Exs_{{Z_1^{n}}}[\Loss(\theta; {Z_1^{n}})]$. In order to estimate this
quantity based
on the data ${Z_1^{\numobs}}$, we solve the convex optimization
problem
%
%e1 #&#
\begin{equation}
\label{EqnMest} {\widehat{\betapar}_{\regpar}}\in\arg\min_{\betapar\in
\real
^\pdim}
\bigl\{ \Loss\bigl(\betapar; {Z_1^{n}}\bigr) + \regpar{
\Regplain(\betapar)} \bigr\},
\end{equation}
where $\regpar> 0$ is a user-defined regularization penalty and
$\Regplain\dvtx{\real^\pdim}\rightarrow\real_+$ is a norm. Note that this
setup allows for the possibility of misspecified models as well.

Our goal is to provide general techniques for deriving bounds on the
difference between any solution ${\widehat{\betapar}_{\regpar}}$ to
the convex
program (\ref{EqnMest}) and the unknown vector $\theta^{*}$. In this
paper we derive bounds on the quantity ${\|{\widehat{\betapar
}_{\regpar}}- {\betapar^*}\|
}$, where the error norm ${\|\cdot\|}$ is induced by some
inner product $\langle{\cdot}, {\cdot} \rangle$ on $\real^\pdim$.
Most often,
this error norm will either be the Euclidean $\ell_2$-norm on vectors
or the analogous Frobenius norm for matrices, but our theory also
applies to certain types of weighted norms. In addition, we provide
bounds on the quantity ${\Regplain({\widehat{\betapar}_{\regpar}}-
{\betapar^*})}$, which measures
the error in the regularizer norm. In the classical setting, the
ambient dimension $\pdim$ stays fixed while the number of observations
$\numobs$ tends to infinity. Under these conditions, there are
standard techniques for proving consistency and asymptotic normality
for the error ${\widehat{\betapar}_{\regpar}}- {\betapar^*}$. In
contrast, the analysis of
this paper is all within a high-dimensional framework, in which the
tuple $(\numobs, \pdim)$, as well as other problem parameters, such as
vector sparsity or matrix rank, etc., are all allowed to tend to
infinity. In contrast to asymptotic statements, our goal is to obtain
explicit finite sample error bounds that hold with high probability.

%s2.2 #&#
\subsection{Decomposability of $\Regplain$}
\label{DecompSection}

The first ingredient in our analysis is a property of the regularizer
known as decomposability, defined in terms of a pair of subspaces
$\ModelSet\subseteq\overline{\ModelSet}$ of $\real^\pdim$. The
role of the
\textit{model subspace} $\ModelSet$ is to capture the constraints
specified by the model; for instance, it might be the subspace of
vectors with a particular support (see Example~\ref{ExaSparseVec}) or
a subspace of low-rank matrices (see Example~\ref{ExaLowRank}). The
orthogonal complement of the space $\overline{\ModelSet}$, namely,
the set
%
%e2 #&#
\begin{equation}\hspace*{8pt}
\overline{\ModelSet}{}^\perp:=\bigl\{ v \in\real^{\pdim} \mid
\langle{u}, {v} \rangle= 0 \mbox{ for all $u \in\overline{\ModelSet}$}
\bigr\},
\end{equation}
is referred to as the \textit{perturbation subspace}, representing
deviations away from the model subspace $\ModelSet$. In the ideal
case, we have \textit{$\overline{\ModelSet}{}^\perp= \ModelSet^\perp
$}, but our definition
allows for the possibility that $\overline{\ModelSet}$ is strictly
larger than
$\ModelSet$, so that $\overline{\ModelSet}{}^\perp$ is strictly
smaller than
${\ModelSet^\perp}$. This generality is needed for treating the case of
low-rank matrices and nuclear norm, as discussed in
Example~\ref{ExaLowRank} to follow.
%
%de1 #&#
\begin{defns}
Given a pair of subspaces $\ModelSet\subseteq\overline{\ModelSet}$,
a norm-based
regularizer $\Regplain$ is \textit{decomposable} with respect to
$(\ModelSet, \overline{\ModelSet}{}^\perp)$ if
%
%e3 #&#
\begin{eqnarray}
\label{EqnDecompose} {\Regplain(\TmpVarA+ \TmpVarB)} &=& {\Regplain
(\TmpVarA)} + {
\Regplain(\TmpVarB)} \nonumber\\[-8pt]\\[-8pt]
&&\eqntext{\mbox{for all $\TmpVarA\in\ModelSet$ and $\TmpVarB
\in
\overline{\ModelSet}{}^\perp$.}}
\end{eqnarray}
\end{defns}

In order to build some intuition, let us consider the ideal
case $\ModelSet= \overline{\ModelSet}$ for the time being, so that the
decomposition (\ref{EqnDecompose}) holds for all pairs $(\TmpVarA,
\TmpVarB) \in\ModelSet\times{\ModelSet^\perp}$. For any given pair
$(\TmpVarA, \TmpVarB)$ of this form, the vector $\TmpVarA+ \TmpVarB$
can be interpreted as a perturbation of the model vector $\TmpVarA$ away
from the subspace $\ModelSet$, and it is desirable that the
regularizer penalize such deviations as much as possible. By the
triangle inequality for a norm, we always have ${\Regplain(\TmpVarA+
\TmpVarB)} \leq{\Regplain(\TmpVarA)} + {\Regplain(\TmpVarB)}$, so
that the
decomposability condition (\ref{EqnDecompose}) holds if and only if the
triangle inequality is tight for all pairs $(\TmpVarA, \TmpVarB) \in
(\ModelSet, {\ModelSet^\perp})$. It is exactly in this setting that the
regularizer penalizes deviations away from the model subspace
$\ModelSet$ as much as possible.\looseness=1

In general, it is not difficult to find subspace pairs that satisfy
the decomposability property. As a trivial example, any regularizer
is decomposable with respect to $\ModelSet= \real^\pdim$ and its
orthogonal complement ${\ModelSet^\perp}= \{0\}$. As will be clear in
our main theorem, it is of more interest to find subspace pairs in
which the model subspace $\ModelSet$ is ``small,'' so that the
orthogonal complement $\ModelSet^\perp$ is ``large.'' To formalize
this intuition, let us define the projection operator
%
%e4 #&#
\begin{equation}
\label{EqnDefnProj} {\Pi_{\ModelSet}(u)} :=\arg\min_{v \in\ModelSet}
\|u - v\|
\end{equation}
with the projection $\Pi_{{\ModelSet^\perp}}$ defined in an analogous
manner. To simplify notation, we frequently use the shorthand
$u_{\ModelSet} = {\Pi_{\ModelSet}(u)}$ and $u_{\ModelSet^\perp} =
\Pi_{{\ModelSet^\perp}}(u)$.

Of interest to us are the action of these projection operators on the
unknown parameter $\param^*\in\real^\pdim$. In the most
desirable setting, the model subspace $\ModelSet$ can be chosen such
that $\param^*_{\ModelSet} \approx\param^*$ or, equivalently,
such that $\param^*_{\ModelSet^\perp} \approx0$. If this can be
achieved with the model subspace $\ModelSet$ remaining relatively
small, then our main theorem guarantees that it is possible to
estimate $\param^*$ at a relatively fast rate. The following
examples illustrate suitable choices of the spaces $\ModelSet$ and
$\overline{\ModelSet}$ in three concrete settings, beginning with the
case of sparse
vectors.
%
%ex1 #&#
\begin{exas}[(Sparse vectors and $\ell_1$-norm regularization)]
\label{ExaSparseVec}
Suppose the
error norm ${\|\cdot\|}$ is the usual $\ell_2$-norm and that the
model class of interest is the set of $\kdim$-sparse vectors in
$\pdim$ dimensions. For any particular subset $\Sset\subseteq
\{1, 2, \ldots, \pdim\}$ with cardinality $\kdim$, we define the
model subspace
%
%e5 #&#
\begin{equation}
\label{EqnModelSparse} \ModelSet(\Sset) :=\bigl\{\TmpVarA\in
\real^\pdim\mid\TmpVarA_j = 0 \mbox{ for all $j \notin S$}
\bigr\}.
\end{equation}
Here our notation reflects the fact that $\ModelSet$ depends
explicitly on the chosen subset $\Sset$. By construction, we have
$\Pi_{\ModelSet(\Sset)}(\param^*) = \param^*$ for any vector
$\param^*$ that is supported on $\Sset$.

In this case, we may define $\overline{\ModelSet}(\Sset) = \ModelSet
(\Sset)$ and
note that the orthogonal complement with respect to the Euclidean
inner product is given by
%
%e6 #&#
\begin{eqnarray}
\label{EqnBsetSparse} \overline{\ModelSet}{}^\perp(\Sset) &=&
\ModelSet^\perp(\Sset)\nonumber\\[-8pt]\\[-8pt]
&=& \bigl\{ \TmpVarB\in\real^\pdim\mid
\TmpVarB_j = 0 \mbox{ for all $j \in S$} \bigr\}.\nonumber
\end{eqnarray}
This set corresponds to the perturbation subspace, capturing
deviations away from the set of vectors with support $\Sset$. We claim
that for any subset $\Sset$, the $\ell_1$-norm ${\Regplain(\TmpVarA
)} = \|
\TmpVarA\|_1$ is decomposable with respect to the pair
$(\ModelSet(\Sset), {\ModelSet^\perp}(\Sset))$. Indeed, by construction
of the subspaces, any $\TmpVarA\in\ModelSet(\Sset)$ can be written
in the partitioned form $\TmpVarA= (\TmpVarA_S, 0_{{{\Sset^c}}})$, where
$\TmpVarA_S \in\real^{\kdim}$ and $0_{{{\Sset^c}}} \in\real^{\pdim-
\kdim}$ is a vector of zeros. Similarly, any vector $\TmpVarB\in
{\ModelSet^\perp}(\Sset)$ has the partitioned representation
$(0_\Sset,
\TmpVarB_{{\Sset^c}})$. Putting together the pieces, we obtain
\[
\| \TmpVarA+ \TmpVarB\|_1 = \bigl\| (\TmpVarA_\Sset, 0) + (0,
\TmpVarB_{{\Sset^c}}) \bigr\|_1 = \| \TmpVarA\|_1 + \|
\TmpVarB\|_1,
\]
showing that the $\ell_1$-norm is decomposable as claimed.
\end{exas}

As a follow-up to the previous example, it is also worth noting that
the same argument shows that for a strictly positive weight vector
$\omega$, the \textit{weighted $\ell_1$-norm} $\|\TmpVarA\|_\omega
:=\sum_{j=1}^\pdim\omega_j |\TmpVarA_j|$ is also decomposable
with respect to the pair $(\ModelSet(\Sset), \overline{\ModelSet
}(\Sset))$. For
another natural extension, we now turn to the case of sparsity models
with more structure.
%
%ex2 #&#
\begin{exas}[(Group-structured norms)]
\label{ExaBlockSparse}
In many applications sparsity arises in a more structured fashion,
with groups of coefficients likely to be zero (or nonzero)
simultaneously. In order to model this behavior, suppose that the
index set $\{1, 2, \ldots, \pdim\}$ can be partitioned into a set of
$\numgroup$ disjoint groups, say, $\GROUP= \{{G}_1, {G}_2,
\ldots, {G}_\numgroup\}$. With this setup, for a given vector
$\SUPERGVEC= (\gqpar_1, \ldots, \gqpar_\numgroup) \in[1,
\infty]^\numgroup$, the associated \textit{$(1,\SUPERGVEC)$-group norm}
takes the form
%
%e7 #&#
\begin{equation}
\label{EqnDefnQblockNorm} {\|\TmpVarA\|_{\GROUP, \SUPERGVEC}} :=\sum
_{t=1}^\numgroup\|\TmpVarA_{{G}_t}
\|_{\gqpar_t}.
\end{equation}
For instance, with the choice $\SUPERGVEC= (2, 2, \ldots, 2)$, we
obtain the group $\ell_1/\ell_2$-norm, corresponding to the
regularizer that underlies the group Lasso \cite{YuaLi06}. On the
other hand, the choice $\SUPERGVEC= (\infty, \ldots, \infty)$,
corresponding to a form of block $\ell_1/\ell_\infty$-regularization,
has also been studied in past work \cite{Turlach05,NegWai08,ZhaRoc06}.
Note that for $\SUPERGVEC= (1, 1, \ldots, 1)$, we obtain the standard
$\ell_1$-penalty. Interestingly, our analysis shows that setting
$\SUPERGVEC\in[2,\infty]^\numgroup$ can often lead to superior
statistical performance.

We now show that the norm ${\|\cdot\|_{\GROUP, \SUPERGVEC}}$ is again
decomposable with respect to appropriately defined subspaces. Indeed,
given any subset $\GroupSset\subseteq\{1, \ldots,
\numgroup\}$ of group indices, say, with cardinality
$\groupcard= |\GroupSset|$, we can define the subspace
%
%e8 #&#
\begin{equation}
\ModelSet(\GroupSset) :=\bigl\{ \TmpVarA\in\real^{\pdim} \mid
\TmpVarA_{{G}_t} = 0 \mbox{ for all $t \notin\GroupSset$} \bigr\}
\end{equation}
as well as its orthogonal complement with respect to the usual
Euclidean inner product
%
%e9 #&#
\begin{eqnarray}
\ModelSet^\perp(\GroupSset) &=& \overline{\ModelSet}{}^\perp(
\GroupSset) \nonumber\\[-8pt]\\[-8pt]
:\!&=&\bigl\{ \TmpVarA\in\real^{\pdim} \mid
\TmpVarA_{{G}_t} = 0 \mbox{ for all $t \in\GroupSset$} \bigr\}.\nonumber
\end{eqnarray}
With these definitions, for any pair of vectors $\TmpVarA\in
\ModelSet(\GroupSset)$ and $\TmpVarB\in\overline{\ModelSet}{}^\perp
(\GroupSset)$, we
have
%
%e10 #&#
\begin{eqnarray}
{\|\TmpVarA+ \TmpVarB\|_{\GROUP, \SUPERGVEC}} &=& \sum_{t \in
\GroupSset}
\|\TmpVarA_{{G}_t} + 0_{{G}_t} \|_{\gqpar_t} \nonumber\\
&&{}+ \sum
_{t \notin\GroupSset} \|0_{{G}_t} + \TmpVarB_{{G}_t}
\|_{\gqpar_t} \\
&=& {\|\TmpVarA\|_{\GROUP, \SUPERGVEC}} + {\|\TmpVarB\|
_{\GROUP, \SUPERGVEC}},\nonumber
\end{eqnarray}
thus verifying the decomposability condition. %
\end{exas}

In the preceding example, we exploited the fact that the groups were
nonoverlapping in order to establish the decomposability property.
Therefore, some modifications would be required in order to choose the
subspaces appropriately for overlapping group regularizers proposed in
past work \cite{Jac09,JenMai10}.
%
%ex3 #&#
\begin{exas}[(Low-rank matrices and nuclear norm)]
\label{ExaLowRank}
Now suppose that each
parameter $\Theta\in\real^{\pdima\times\pdimb}$ is a matrix; this
corresponds to an instance of our general setup with $\pdim= \pdima
\pdimb$, as long as we identify the space $\real^{\pdima\times
\pdimb}$ with $\real^{\pdima\pdimb}$ in the usual way. We equip
this space with the inner product ${\llangle{\Theta}, {\Gamma
}\rrangle} :=
\trace(\Theta\Gamma^T)$, a choice which yields (as the induced norm)
the \textit{Frobenius norm}
%
%e11 #&#
\begin{equation}
{\norm{\Theta}\norm_{\mathrm{F}}} :=\sqrt{{\llangle{\Theta}, {\Theta}
\rrangle}} = \sqrt{\sum_{j=1}^\pdima
\sum_{k=1}^\pdimb\Theta_{jk}^2}.
\end{equation}
In many settings, it is natural to consider estimating matrices that
are low-rank; examples include principal component analysis, spectral
clustering, collaborative filtering and matrix completion. With
certain exceptions, it is computationally expensive to enforce a
rank-constraint in a direct manner, so that a variety of researchers
have studied the \textit{nuclear norm}, also known as the trace norm, as
a surrogate for a rank constraint. More precisely, the nuclear norm
is given by
%
%e12 #&#
\begin{equation}
{\norm{\Theta}\norm}_{\mathrm{nuc}} :=\sum_{j=1}^{\min\{
\pdima, \pdimb\}}
\sigma_j(\Theta),
\end{equation}
where $\{\sigma_j(\Theta) \}$ are the singular values of the matrix
$\Theta$.

The nuclear norm is decomposable with respect to appropriately chosen
subspaces. Let us consider the class of matrices $\Theta\in
\real^{\pdima\times\pdimb}$ that have rank $r \leq\min
\{\pdima, \pdimb\}$. For any given matrix $\Theta$, we let
$\rowsp(\Theta) \subseteq\real^\pdimb$ and $\colsp(\Theta)
\subseteq
\real^\pdima$ denote its row space and column space, respectively. Let
$U$ and $V$ be a given pair of $r$-dimensional subspaces $U \subseteq
\real^\pdima$ and $V \subseteq\real^\pdimb$; these subspaces will
represent left and right singular vectors of the target matrix
$\Theta^*$ to be estimated. For a given pair $(U, V)$, we can
define\vspace*{1pt}
the subspaces $\ModelSet(U, V)$ and $\overline{\ModelSet}{}^\perp(U,
V)$ of
$\real^{\pdima\times\pdimb}$ given by
%
%e13 #&#
\begin{subequation}
%
%e13.a #&#
%e13.b #&#
\begin{eqnarray}\hspace*{12pt}
\ModelSet(U, V)
&:=&\bigl\{ \Theta\in\real^{\pdima\times
\pdimb} \mid\rowsp(
\Theta) \subseteq V,\nonumber\\[-8pt]\\[-8pt]
&&\hspace*{71.3pt} \colsp(\Theta) \subseteq U \bigr\}\nonumber
\end{eqnarray}
and
\begin{eqnarray}
\overline{\ModelSet}{}^\perp(U, V) &:=& \bigl\{ \Theta\in
\real^{\pdima\times
\pdimb} \mid\rowsp(\Theta) \subseteq V^\perp,\hspace*{-32pt}\nonumber\\[-8pt]\\[-8pt]
&&\hspace*{72pt} \colsp(
\Theta) \subseteq U^\perp\bigr\}.\hspace*{-32pt}\nonumber
\end{eqnarray}
\end{subequation}
So as to simplify notation, we omit the indices $(U, V)$ when they are
clear from context. Unlike the preceding examples, in this case, the
set $\ModelSet$ is not\footnote{However, as is required by our theory,
we do have the inclusion $\ModelSet\subseteq\overline{\ModelSet}$.
Indeed, given any $\Theta\in\ModelSet$ and
$\Gamma\in\overline{\ModelSet }{}^\perp$, we have $\Theta^T \Gamma= 0$
by definition, which implies that ${\llangle{\Theta}, {\Gamma}\rrangle}
= \trace(\Theta^T \Gamma) = 0$. Since
$\Gamma\in\overline{\ModelSet}{}^\perp$ was arbitrary, we have shown
that $\Theta$ is orthogonal to the space
$\overline{\ModelSet}{}^\perp$, meaning that it must belong
to~$\overline{\ModelSet}$.} equal to $\overline{\ModelSet}$.

Finally, we claim that the nuclear norm is decomposable with respect
to the pair $(\ModelSet, \overline{\ModelSet}{}^\perp)$. By
construction, any pair of
matrices $\Theta\in\ModelSet$ and $\Gamma\in\overline{\ModelSet
}{}^\perp$ have
orthogonal row and column spaces, which implies the required
decomposability condition---name\-ly, $\norm\Theta+\Gamma\norm_{{1}} =
\norm\Theta\norm_{{1}} + \norm\Gamma\norm_{{1}}$.
\end{exas}

A line of recent work
(e.g., \cite{Chand09,XuCaSa2010,CanLi10,AgaNegWai11,HUPAPER,MccoyTr2010})
has studied matrix problems\vadjust{\eject} involving the sum of a low-rank matrix
with a sparse matrix, along with the regularizer formed by a weighted
sum of the nuclear norm and the elementwise $\ell_1$-norm. By a
combination of Examples~\ref{ExaSparseVec} and
\ref{ExaLowRank}, this regularizer also satisfies the
decomposability property with respect to appropriately defined
subspaces.

%s2.3 #&#
\subsection{A Key Consequence of Decomposability}

Thus far, we have specified a class (\ref{EqnMest}) of $M$-esti\-mators
based on regularization, defined the notion of decomposability for the
regularizer and worked through several illustrative examples. We now
turn to the statistical consequences of decomposability---more
specifically, its implications for the error vector $\delparhatlam=
{\widehat{\betapar}_{\regpar}}- {\betapar^*}$, where $\widehat
{\betapar}\in\real^\pdim$ is any
solution of the regularized $M$-estimation procedure (\ref{EqnMest}).
For a given inner product $\langle{\cdot}, {\cdot} \rangle$, the
dual norm of
$\Regplain$ is given by
%
%e14 #&#
\begin{equation}
{\Regplain^*(v)} :=\sup_{u \in\real^\pdim\setminus\{0\}} \frac
{\langle{u}, {v} \rangle}{{\Regplain(u)}} =
\sup_{{\Regplain(u)} \leq1} \langle{u}, {v} \rangle.
\end{equation}
This notion is best understood by working through some examples.

%pa2.3.subsubsection.1 #&#
\subsubsection*{Dual of $\ell_1$-norm} For the $\ell_1$-norm
${\Regplain(u)} = \|u\|_1$ previously discussed in
Example~\ref{ExaSparseVec}, let us compute its dual norm with respect
to the Euclidean inner product on $\real^\pdim$. For any vector $v
\in
\real^{\pdim}$, we have
\begin{eqnarray*}
\sup_{\|u\|_1 \leq1} \langle{u}, {v} \rangle&\leq&\sup_{\|u\|
_1 \leq1} \sum
_{k=1}^\pdim|u_k|
|v_k| \\
&\leq&\sup_{\|u\|_1 \leq1} \Biggl( \sum
_{k=1}^\pdim|u_k| \Biggr)
\max_{k = 1, \ldots, \pdim} |v_k| \\
&=& \|v\|_\infty.
\end{eqnarray*}
We claim that this upper bound actually holds with equality. In
particular, letting $j$ be any index\break for~which $|v_j|$ achieves the
maximum $\|v\|_\infty= \break\max_{k = 1, \ldots, \pdim} |v_k|$, suppose
that we form a vector ${\overline{u}}\in\real^\pdim$ with
${\overline{u}}_j =
\sign(v_j)$ and ${\overline{u}}_k = 0$ for all $k \neq j$. With this
choice, we have $\|{\overline{u}}\|_1 \leq1$ and, hence,
\[
\sup_{\|u\|_1 \leq1} \langle{u}, {v} \rangle\geq\sum
_{k=1}^\pdim{\overline{u}}_k
v_k = \|v\|_\infty,
\]
showing that the dual of the $\ell_1$-norm is the $\ell_\infty$-norm.

%
%f1 #&#
\begin{figure*}

\includegraphics{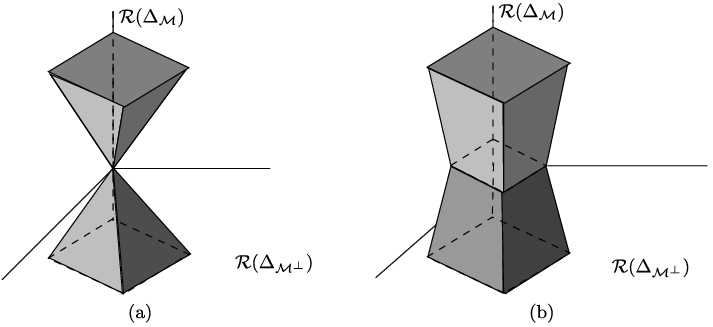}

\caption{Illustration of the set $\ConvSet(\ModelSet, {\ModelSet^\perp};
{\betapar^*})$ in the special case $\delpar= (\delpar_1, \delpar_2,
\delpar_3) \in\real^3$ and regularizer ${\Regplain(\delpar)} =
\|\delpar\|_1$, relevant for sparse vectors
(Example \protect\ref{ExaSparseVec}). This picture shows the case
$\Sset=
\{3\}$, so that the model subspace is $\ModelSet(\Sset) = \{\delpar
\in\real^3 \mid\delpar_1 = \delpar_2 = 0 \}$ and its
orthogonal complement is given by ${\ModelSet^\perp}(\Sset) =
\{\delpar\in\real^3 \mid\delpar_3 = 0 \}$. \textup{(a)} In the
special case when $\param^*_1 = \param^*_2 = 0$, so that
$\param^*\in\ModelSet$, the set $\ConeSet(\ModelSet,
{\ModelSet^\perp}; \param^*)$ is a cone. \textup{(b)} When $\param^*$ does
not belong to $\ModelSet$, the set $\ConvSet(\ModelSet,
{\ModelSet^\perp}; {\betapar^*})$ is enlarged in the coordinates
$(\delpar_1, \delpar_2)$ that span~${\ModelSet^\perp}$. It is no longer
a cone, but is still a star-shaped set.}
\label{FigConeSet}
\end{figure*}
%

%pa2.3.subsubsection.2 #&#
\subsubsection*{Dual of group norm} Now recall the group norm from
Example~\ref{ExaBlockSparse}, specified in terms of a vector
$\SUPERGVEC\in[2, \infty]^{\numgroup}$. A similar calculation shows
that its dual norm, again with respect to the Euclidean norm on
$\real^\pdim$, is given by
%
%e15 #&#
\begin{eqnarray}
{\|v\|_{\GROUP, \SUPERGVEC^*}} &=& \max_{t = 1, \ldots, \numgroup} \|v\|
_{\gqpar^*_t} \hspace*{35pt}\nonumber\\[-8pt]\\[-8pt]
&&\eqntext{\mbox{where $\displaystyle
\frac{1}{\gqpar_t} + \frac{1}{\gqpar^*_t} = 1 $ are dual exponents.}}
\end{eqnarray}
As special cases of this general duality relation, the block $(1,2)$
norm that underlies the usual group Lasso leads to a block $(\infty,
2)$ norm as the dual, whereas the block $(1, \infty)$ norm leads
to a block $(\infty, 1)$ norm as the dual.
%

%pa2.3.subsubsection.3 #&#
\subsubsection*{Dual of nuclear norm} For the nuclear norm, the
dual is defined with respect to the trace inner product on the space
of matrices. For any matrix $N \in\real^{\pdima\times\pdimb}$, it
can be shown that
\begin{eqnarray*}
{\Regplain^*(N)} &=& \sup_{{\norm{M}\norm}_{\mathrm{nuc}} \leq1} {\llangle
{M}, {N}\rrangle} = \norm{N}
\norm_{\mathrm{op}} \\
&=& \max_{j = 1, \ldots,
\min\{\pdima, \pdimb\}} \sigma_j(N),
\end{eqnarray*}
corresponding to the $\ell_\infty$-norm applied to the vector
$\sigma(N)$ of singular values. In the special case of diagonal
matrices, this fact reduces to the dual relationship between the
vector $\ell_1$ and $\ell_\infty$-norms.

The dual norm plays a key role in our general theory, in particular, by
specifying a suitable choice of the regularization weight $\regpar$.
We summarize in the following:
%
%le1 #&#
\begin{lems}
\label{LemOptBound}
Suppose that $\Loss$ is a convex and differentiable function, and
consider any optimal solution $\widehat{\betapar}$ to the optimization
problem (\ref{EqnMest}) with a strictly positive regularization
parameter satisfying
%
%e16 #&#
\begin{equation}
\label{EqnLamBound} \regpar\geq2 {\Regplain^*\bigl(\nabla\Loss\bigl({
\betapar^*}; {Z_1^{n}}\bigr)\bigr)}.
\end{equation}
Then\vspace*{1pt} for any pair $(\ModelSet, \overline{\ModelSet}{}^\perp)$ over
which $\Regplain$ is
decomposable, the error ${\widehat{\delpar}}= {\widehat{\betapar
}_{\regpar}}- {\betapar^*}$ belongs
to the set
%
%e17 #&#
\begin{eqnarray}
\label{EqnDefnConeSet}
&&\ConvSet\bigl(\ModelSet, \overline{\ModelSet
}{}^\perp;
{\betapar^*}\bigr) \nonumber\\
&&\quad:=\bigl\{ \delpar\in\real^\pdim\mid{\Regplain(
\delpar_{\bar{\ModelSet}^\perp})} \\
&&\quad\hspace*{18.1pt}\leq3 {\Regplain(\delpar_{\bar
{\ModelSet}})} + 4 {\Regplain
\bigl(\param^*_{\ModelSet^\perp}\bigr)} \bigr\}.\nonumber
\end{eqnarray}
\end{lems}

We prove this result in the supplementary
appen\-dix~\cite{Neg12supp}. It has the following important
consequence: for any decomposable regularizer and an appropriate
choice (\ref{EqnLamBound}) of regularization\vspace*{1pt} parameter, we are
guaranteed that the error vector ${\widehat{\delpar}}$ belongs to a very
specific set, depending on the unknown vector ${\betapar^*}$. As
\mbox{illustrated} in Figure~\ref{FigConeSet}, the geometry of the set
$\ConvSet$ depends on the relation between ${\betapar^*}$ and the model
subspace $\ModelSet$. When ${\betapar^*}\in\ModelSet$, then we are
guaranteed that ${\Regplain(\param^*_{\ModelSet^\perp})} = 0$. In
this case,\vspace*{1pt}
the constraint (\ref{EqnDefnConeSet}) reduces to
${\Regplain(\delpar_{\bar{\ModelSet}^\perp})} \leq3 {\Regplain
(\delpar_{\bar{\ModelSet}})}$, so
that $\ConvSet$ is a cone, as illustrated in panel (a). In the more
general case when ${\betapar^*}\notin\ModelSet$ so that
${\Regplain(\param^*_{\ModelSet^\perp})} \neq0$, the set
$\ConvSet
$ is
\textit{not} a cone, but rather a star-shaped set [panel (b)]. As will
be clarified in the sequel, the case ${\betapar^*}\notin\ModelSet$
requires a more delicate treatment.

%%%%%%%%%%%%%%%%%%%%%%%%%%%%%%%%%%%%%%%%%%%%%%%%%%%%%%%%%%%%%%%%%%%%
%s2.4 #&#
\subsection{Restricted Strong Convexity}
\label{RSCSec}

\begin{figure*}

\includegraphics{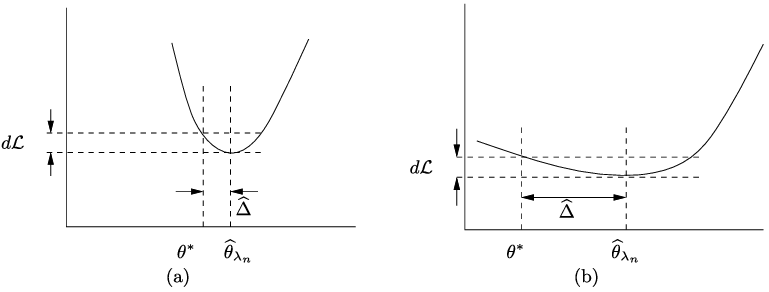}

\caption{Role of curvature in distinguishing parameters. \textup{(a)} Loss
function has high curvature around ${\widehat{\delpar}}$. A small excess
loss $d \Loss= |\Loss({\widehat{\betapar}_{\regpar}}) - \Loss
({\betapar^*})|$ guarantees
that the parameter error $\widehat{\Delta}= {\widehat{\betapar
}_{\regpar}}- {\betapar^*}$ is
also small. \textup{(b)} A less desirable setting, in which the loss
function has relatively low curvature around the optimum.}
\label{FigCurvature}\vspace*{-3pt}
\end{figure*}

We now turn to an important requirement of the loss function and its
interaction with the statistical model. Recall that $\widehat{\Delta}=
{\widehat{\betapar}_{\regpar}}- {\betapar^*}$ is the difference
between an optimal solution
${\widehat{\betapar}_{\regpar}}$ and the true parameter, and
consider the loss
difference\footnote{To simplify notation, we frequently write
$\Loss(\betapar)$ as shorthand for $\Loss(\betapar; {Z_1^{n}})$ when
the underlying data ${Z_1^{n}}$ is clear from context.}
$\Loss({\widehat{\betapar}_{\regpar}})- \Loss({\betapar^*})$. In
the classical setting,
under fairly mild conditions, one expects that the loss
difference should converge to zero as the sample size $\numobs$
increases. It is important to note, however, that such convergence on
its own is \textit{not sufficient} to guarantee that ${\widehat
{\betapar}_{\regpar}}$ and
${\betapar^*}$ are close or, equivalently, that $\widehat{\Delta}$
is small.
Rather, the closeness depends on the curvature of the loss function,
as illustrated in Figure~\ref{FigCurvature}.
%
%f2 #&#
In a desirable setting [panel (a)], the loss function is sharply
curved around its optimum ${\widehat{\betapar}_{\regpar}}$, so that
having a small loss
difference $|\Loss({\betapar^*}) - \Loss({\widehat{\betapar
}_{\regpar}})|$ translates to a
small error $\widehat{\Delta}= {\widehat{\betapar}_{\regpar}}-
{\betapar^*}$. Panel (b)
illustrates a less desirable setting, in which the loss function is
relatively flat, so that the loss difference can be small while the
error $\widehat{\Delta}$ is relatively large.

The standard way to ensure that a function is ``not too flat'' is via
the notion of strong convexity. Since $\Loss$ is differentiable by
assumption, we may perform a first-order Taylor series expansion at
${\betapar^*}$ and in some direction $\Delta$; the error in this Taylor
series is given~by
%
%e18 #&#
\begin{eqnarray}
\label{EqnStrongConvexity}
\delLike\bigl(\delpar, \theta^{*}\bigr)
&:=&\Loss\bigl({\betapar^*}+ \delpar\bigr) - \Loss\bigl({\betapar
^*}\bigr)
\nonumber\\[-8pt]\\[-8pt]
&&{}
-
\bigl\langle{\nabla\Loss\bigl({\betapar^*}\bigr)}, {\delpar} \bigr
\rangle.\nonumber%\hspace*{-22pt}%\nonumber
\end{eqnarray}
One way in which to enforce that $\Loss$ is strongly convex is to
require the existence of some positive constant $\curve> 0$ such that
$\delLike(\delpar, \theta^{*}) \geq\curve{\|\delpar\|}^2$ for all
$\delpar\in\real^\pdim$ in a neighborhood of ${\betapar^*}$. When the
loss function is twice differentiable, strong convexity amounts to
lower bound on the eigenvalues of the Hessian $\nabla^2
\Loss(\theta)$, holding uniformly for all $\theta$ in a neighborhood
of $\param^*$.

Under classical ``fixed $\pdim$, large $\numobs$'' scaling, the loss
function will be strongly convex under mild conditions. For instance,
suppose that population risk ${\overline{\Loss}}$ is strongly convex or,
equivalently, that the Hessian $\nabla^2 {\overline{\Loss}}(\theta
)$ is strictly
positive definite in a neighborhood of $\param^*$. As a concrete
example, when the loss function $\Loss$ is defined based on negative
log likelihood of a statistical model, then the Hessian $\nabla^2
{\overline{\Loss}}(\theta)$ corresponds to the Fisher information
matrix, a
quantity which arises naturally in asymptotic statistics. If the
dimension $\pdim$ is fixed while the sample size $\numobs$ goes to
infinity, standard arguments can be used to show that (under mild
regularity conditions) the random Hessian $\nabla^2 \Loss(\theta)$
converges to $\nabla^2 {\overline{\Loss}}(\theta)$ uniformly for
all $\theta$ in
an open neighborhood of~$\param^*$. In contrast, whenever the pair
$(\numobs, \pdim)$ both increase in such a way that $\pdim> \numobs$,
the situation is drastically different: the Hessian matrix $\nabla^2
\Loss(\theta)$ is often singular. As a concrete example, consider
linear regression based on samples $Z_i = (y_i, x_i) \in\real
\times\real^\pdim$, for $i = 1, 2, \ldots, \numobs$. Using the
least squares loss $\Loss(\theta) = \frac{1}{2 \numobs} \|y - X
\theta\|_2^2$, the $\pdim\times\pdim$ Hessian matrix $\nabla^2
\Loss(\theta) = \frac{1}{\numobs} X^T X$ has rank at most
$\numobs$, meaning that the loss cannot be strongly convex when $\pdim
> \numobs$. Consequently, it impossible to guarantee global strong
convexity, so that we need to restrict the set of directions $\Delta$
in which we require a curvature condition.

%f3 #&#
\begin{figure*}

\includegraphics{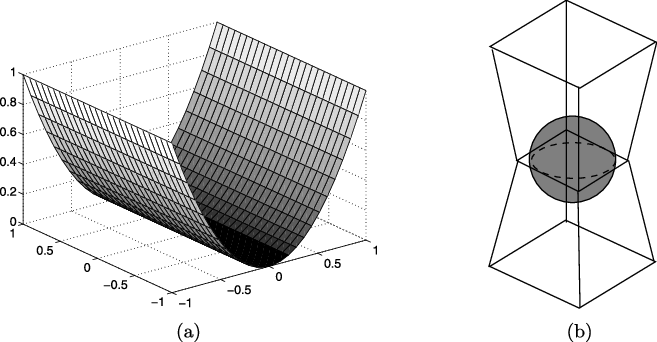}

\caption{\textup{(a)} Illustration of a generic loss function in the
high-dimensional $\pdim> \numobs$ setting: it is curved in certain
directions, but completely flat in others. \textup{(b)} When $\param^*
\notin\ModelSet$, the set $\ConvSet(\ModelSet, \overline{\ModelSet
}{}^\perp;
\param^*)$
contains a ball centered at the origin, which necessitates a
tolerance term $\rscslop({\betapar^*}) > 0$ in the definition of
restricted strong convexity.}
\label{FigRSC}
\end{figure*}

Ultimately, the only direction of interest is given by the error
vector $\widehat{\delpar}= {\widehat{\betapar}_{\regpar}}-
{\betapar^*}$. Recall that
Lemma~\ref{LemOptBound} guarantees that, for suitable choices of the
regularization parameter $\regpar$, this error vector must belong to
the set $\ConvSet(\ModelSet, \overline{\ModelSet}{}^\perp; {\betapar
^*})$, as previously
defined (\ref{EqnDefnConeSet}). Consequently, it suffices to ensure
that the function is strongly convex over this set, as formalized
in the following:
%
%de2 #&#
\begin{defns}
\label{DefnRSC}
The loss function satisfies a \textit{restricted strong
convexity} (RSC) condition with \textit{curvature $\NewMcon>
0$} and \textit{tolerance function $\rscslop$} if
%
%e19 #&#
\begin{eqnarray}
\label{EqnDefnRSC}
\delLike\bigl(\delpar, \theta^{*}\bigr) &\geq&
\NewMcon{\|{\delpar}\|^2} - \rscslop^2\bigl(
\theta^{*}\bigr)\nonumber\\[-8pt]\\[-8pt]
&&\eqntext{\mbox{for all $\delpar\in\ConvSet\bigl(\ModelSet,
\overline{\ModelSet}{}^\perp; {\betapar^*}\bigr)$.}}
\end{eqnarray}
\end{defns}

In the simplest of cases---in particular, when $\theta^{*}\in
\ModelSet$---there are many statistical models for which this RSC
condition holds with tolerance $\rscslop(\param^*) = 0$. In the
more general setting, it can hold only with a nonzero tolerance term,
as illustrated in Figure~\ref{FigRSC}(b). As our proofs will clarify,
we in fact require only the lower bound (\ref{EqnDefnRSC}) to hold for
the intersection of $\ConvSet$ with a local ball $\{ \|\Delta\| \leq
R \}$ of some radius centered at zero. As will be clarified later,
this restriction is not necessary for the least squares loss, but is
essential for more general loss functions, such as those that arise in
generalized linear models.

We will see in the sequel that for many loss functions, it is possible
to prove that with high probability the first-order Taylor series
error satisfies a lower bound of the form
%
%e20 #&#
\begin{eqnarray}
\label{EqnSahand}
\delLike\bigl(\delpar, \theta^{*}\bigr) &\geq&
\kappa_1 {\|{\delpar}\|^2} - \GARTWO g(\numobs, \pdim)
\Regplain^2(\Delta) \nonumber\\[-8pt]\\[-8pt]
&&\eqntext{\mbox{for all $\|\Delta\| \leq1$,}}
\end{eqnarray}
where $\kappa_1, \GARTWO$ are positive constants and $g(\numobs,
\pdim)$ is a function of the sample size $\numobs$ and ambient
dimension $\pdim$, decreasing in the sample size. For instance, in
the case of $\ell_1$-regularization, for covariates with suitably
controlled tails, this type of bound holds for the least squares
loss with the function $g(\numobs, \pdim) = \frac{\log
\pdim}{\numobs}$; see equation (\ref{EqnGarvesh}) to follow. For
generalized linear models and the $\ell_1$-norm, a similar type of
bound is given in equation (\ref{EqnSahandGLM}). We also provide a
bound of this form for the least-squares loss group-structured norms
in equation (\ref{EqnREGroupGeneral}), with a different choice of the
function $g$ depending on the group structure.

A bound of the form (\ref{EqnSahand}) implies a form of restricted
strong convexity as long as ${\Regplain(\delpar)}$ is not ``too large''
relative to $\|\delpar\|$. In order to formalize this notion, we
define a quantity that relates the error norm and the regularizer:
%
%de3 #&#
\begin{defns}[(Subspace compatibility constant)]
\label{DefnSubspaceCompat}
For any subspace $\ModelSet$ of $\real^\pdim$, the \textit{subspace
compatibility constant} with respect to the pair $(\Regplain,
{\|\cdot\|})$ is given by
%
%e21 #&#
\begin{equation}
\label{EqnDefnCompat} {\Psi(\ModelSet)} :=\sup_{u \in\ModelSet
\setminus\{0\}}
\frac{{\Regplain(u)}}{{\|u\|}}.
\end{equation}
\end{defns}

This quantity reflects the degree of compatibility between
the regularizer and the error norm over the subspace $\ModelSet$. In
alternative terms, it is the Lipschitz constant of the regularizer
with respect to the error norm, restricted to the subspace
$\ModelSet$. As a simple example, if $\ModelSet$ is a
$\kdim$-dimensional coordinate subspace, with regularizer ${\Regplain
(u)} =
\|u\|_1$ and error norm ${\|u\|} = \|u\|_2$, then we have
${\Psi(\ModelSet)} = \sqrt{\kdim}$.

This compatibility constant appears explicitly in the bounds of our
main theorem and also arises in establishing restricted strong
convexity. Let us now illustrate how it can be used to show that the
condition (\ref{EqnSahand}) implies a form of restricted strong
convexity. To be concrete, let us suppose that $\param^*$ belongs
to a subspace $\ModelSet$; in this case, membership of $\Delta$ in the
set $\ConvSet(\ModelSet, \overline{\ModelSet}{}^\perp; {\betapar^*})$
implies that
${\Regplain(\delpar_{\bar{\ModelSet}^\perp})} \leq3 {\Regplain
(\delpar_{\bar{\ModelSet}})}$.
Consequently, by the triangle inequality and the
definition (\ref{EqnDefnCompat}), we have
\begin{eqnarray*}
{\Regplain(\delpar)} &\leq&{\Regplain(\delpar_{\bar{\ModelSet
}^\perp})} + {\Regplain(
\delpar_{\bar{\ModelSet}})} \leq4 {\Regplain(\delpar_{\bar{\ModelSet
}})} \\
&\leq&
4 {\Psi(\overline{\ModelSet})} \|\delpar\|.
\end{eqnarray*}
Therefore, whenever a bound of the form (\ref{EqnSahand}) holds and
$\param^*\in\ModelSet$, we are guaranteed that
\begin{eqnarray}
\delLike\bigl(\delpar, \theta^{*}\bigr) &\geq&\bigl\{
\kappa_1 - 16 \GARTWO{\Psi^2(\overline{\ModelSet})} g(
\numobs, \pdim) \bigr\} {\|{\delpar}\|^2} \nonumber\\
&&\eqntext{\mbox{for all $\|\delpar\|
\leq1$.}}
\end{eqnarray}
Consequently, as long as the sample size is large enough that $16
\GARTWO{\Psi^2(\overline{\ModelSet})} g(\numobs, \pdim) <
\frac{\kappa_1}{2}$,
the restricted strong convexity condition will hold with $\NewMcon=
\frac{\kappa_1}{2}$ and $\rscslop(\theta^{*}) = 0$. We make use of
arguments of this flavor throughout this paper.

%%%%%%%%%%%%%%%%%%%%%%%%%%%%%%%%%%%%%%%%%%%%%%%%%%%%%%%%%%%%%%%%%%%%%%%%%%%%%%%%%%%%%%%%%%%%%%%%%%%

%s3 #&#
\section{Bounds for General $M$-Estimators}
\label{SecMain}

We are now ready to state a general result that provides bounds and
hence convergence rates for the error ${\|{\widehat{\betapar
}_{\regpar}}- {\betapar^*}\|}$,
where ${\widehat{\betapar}_{\regpar}}$ is any optimal solution of
the convex
program (\ref{EqnMest}). Although it may appear somewhat abstract at
first sight, this result has a number of concrete and useful
consequences for specific models. In particular, we recover as an
immediate corollary the best known results about estimation in sparse
linear models with general designs \cite{BiRiTsy08,MeiYu09}, as well
as a number of new results, including minimax-optimal rates for
estimation under $\ell_q$-sparsity constraints and estimation of
block-structured sparse matrices. In results that we report
elsewhere, we also apply these theorems to establishing results for
sparse generalized linear models~\cite{Neg09}, estimation of low-rank
matrices \mbox{\cite{NegWai10b,NegWai09}}, \mbox{matrix} decomposition
problems \cite{AgaNegWai11} and sparse nonparametric regression
models \cite{RasWaiYu10b}.

Let us recall our running assumptions on the structure of
the convex program (\ref{EqnMest}).
\begin{longlist}[(G2)]
\item[(G1)] The regularizer $\Regplain$ is a norm and is decomposable with
respect to the subspace pair $(\ModelSet, \overline{\ModelSet}{}^\perp
)$, where
$\ModelSet\subseteq\overline{\ModelSet}$.
\item[(G2)] The loss function $\Loss$ is convex and differentiable,
and satisfies restricted strong convexity with curvature $\NewMcon$
and tolerance $\rscslop$.
\end{longlist}
The reader should also recall the definition (\ref{EqnDefnCompat}) of
the subspace compatibility constant. With this\vadjust{\goodbreak} notation, we can now
state the main result of this paper:
%
%th1 #&#
\begin{theos}[(Bounds for general models)]
\label{ThmMain} Under conditions \textup{(G1)} and \textup{(G2)},
consider the problem (\ref{EqnMest}) based on a strictly positive
regularization constant $\regpar\geq2
{\Regplain^*(\nabla\Loss({\betapar^*}))}$. Then any optimal solution
${\widehat{\betapar}_{\regpar}}$ to the convex program (\ref {EqnMest})
satisfies the bound
%
%e22 #&#
\begin{eqnarray}
\label{mainthm}\quad
{\bigl\|{\widehat{\betapar}_{\regpar}}- {\betapar^*}
\bigr\|}^2 &\leq&9 \frac{\regpar^2}{\NewMcon^2} {\Psi^2(\overline{
\ModelSet})} \nonumber\\[-8pt]\\[-8pt]
&&{}+ \frac{\regpar}{\NewMcon} \bigl\{ 2 \rscslop^2\bigl(
\theta^{*}\bigr) + 4 {\Regplain\bigl(\param^*_{\ModelSet^\perp}\bigr)}
\bigr\}.\nonumber
\end{eqnarray}
\end{theos}
\begin{Remarks*}
Let us consider in more detail some
different features of this result.

\begin{longlist}[(a)]
\item[(a)] It should be noted that Theorem~\ref{ThmMain} is actually a
\textit{deterministic} statement about the set of optimizers of the
convex program (\ref{EqnMest}) for a fixed choice of $\regpar$.
Although the program is convex, it need not be strictly convex, so
that the global optimum might be attained at more than one point
${\widehat{\betapar}_{\regpar}}$. The stated bound holds for any of
these optima.
Probabilistic analysis is required when Theorem~\ref{ThmMain} is
applied to particular statistical models, and we need to verify that
the regularizer satisfies the condition
%
%e23 #&#
\begin{equation}
\label{EqnValid} \regpar\geq2 {\Regplain^*\bigl(\nabla\Loss\bigl
({\betapar^*}
\bigr)\bigr)}
\end{equation}
and that the loss satisfies the RSC condition. A challenge here is
that since ${\betapar^*}$ is unknown, it is usually impossible to compute
the right-hand side of the condition (\ref{EqnValid}). Instead, when
we derive consequences of Theorem~\ref{ThmMain} for different
statistical models, we use concentration inequalities in order to
provide bounds that hold with high probability over the data.
\item[(b)] Second, note that Theorem~\ref{ThmMain} actually provides a
\textit{family of bounds}, one for each pair $(\ModelSet, \overline
{\ModelSet}{}^\perp)$
of subspaces for which the regularizer is decomposable. Ignoring
the term involving $\rscslop$ for the moment, for any given pair,
the error bound is the sum of two terms, corresponding to estimation
error $\Esterr$ and approximation error $\Apperr$, given by,
respectively,
%
%e24 #&#
\begin{eqnarray}
\Esterr&:=&9 \frac{\regpar^2}{\NewMcon^2} {\Psi^2(\overline{\ModelSet})}
\quad\mbox{and}\nonumber\\[-8pt]\\[-8pt]
\Apperr&:=&4 \frac{\regpar}{\NewMcon} {\Regplain\bigl(
\param^*_{\ModelSet^\perp}\bigr)}.\nonumber
\end{eqnarray}
As the dimension of the subspace $\ModelSet$ increases (so that the
dimension of $\ModelSet^\perp$ decreases), the approximation error
tends to zero. But since $\ModelSet\subseteq\overline{\ModelSet}$,
the estimation
error is increasing at the same time. Thus, in the usual way, optimal
rates are obtained by choosing $\ModelSet$ and $\overline{\ModelSet
}$ so as to
balance these two contributions to the error. We illustrate such
choices for various specific models to follow.
\item[(c)] As will be clarified in the sequel, many
high-dimensional statistical models have an unidentifiable
component, and the tolerance term $\rscslop$ reflects the degree of
this nonidentifiability.
\end{longlist}

A large body of past work on sparse linear regression has focused on
the case of exactly sparse regression models for which the unknown
regression vector $\param^*$ is $\kdim$-sparse. For this special
case, recall from Example~\ref{ExaSparseVec} in
Section~\ref{DecompSection} that we can define an $\kdim$-dimensional
subspace $\ModelSet$ that contains $\param^*$. Consequently, the
associated set $\ConeSet(\ModelSet, {\ModelSet^\perp}; \param^*)$
is a
cone [see Figure~\ref{FigConeSet}(a)], and it is thus possible to
establish that restricted strong convexity (RSC) holds with tolerance
parameter \mbox{$\rscslop(\theta^{*}) = 0$}. This same reasoning applies to
other statistical models, among them group-sparse regression, in which
a small subset of groups are active, as well as low-rank matrix
estimation. The following corollary provides a simply stated bound
that covers all of these models:\looseness=1
\end{Remarks*}

%
%co1 #&#
\begin{cors}
\label{CorExact}
Suppose that, in addition to the conditions of Theorem~\ref{ThmMain},
the unknown ${\betapar^*}$ belongs to $\ModelSet$ and the RSC condition
holds over $\ConvSet(\ModelSet, \overline{\ModelSet}, \param^*)$ with
$\rscslop(\param^*) = 0$. Then any optimal solution ${\widehat
{\betapar}_{\regpar}}$ to
the convex program (\ref{EqnMest}) satisfies the bounds
%
%e25 #&#
\begin{subequation}
%
%e25.a #&#
%e25.b #&#
\begin{eqnarray}
\label{EqnThreePart} {\bigl\|{\widehat{\betapar}_{\regpar}}- {\betapar^*}\bigr\|}
& \leq & 9 \frac{\regpar^2}{\NewMcon} {\Psi^2(\overline{\ModelSet})}
\end{eqnarray}
and
\begin{eqnarray}
\label{EqnThreePartReg} {\Regplain\bigl({\widehat{
\betapar}_{\regpar}}- {\betapar^*}\bigr)} & \leq & 12 \frac
{\regpar
}{\NewMcon} {
\Psi^2(\overline{\ModelSet})}.
\end{eqnarray}
\end{subequation}
\end{cors}

Focusing first on the bound (\ref{EqnThreePart}), it
consists of three terms, each of which has a natural interpretation.
First, it is inversely proportional to the RSC constant $\NewMcon$, so
that higher curvature guarantees lower error, as is to be expected.
The error bound grows proportionally with the subspace compatibility
constant ${\Psi(\overline{\ModelSet})}$, which measures the
compatibility between
the regularizer $\Regplain$ and error norm ${\|\cdot\|}$ over the
subspace $\overline{\ModelSet}$ (see Definition \ref
{DefnSubspaceCompat}). This term
increases with the size of subspace $\overline{\ModelSet}$, which
contains the model
subspace $\ModelSet$. Third, the bound also scales linearly with the
regularization parameter $\regpar$, which must be strictly positive
and satisfy the lower bound (\ref{EqnValid}). The
bound (\ref{EqnThreePartReg}) on the error measured in the regularizer
norm is similar, except that it scales quadratically with the subspace
compatibility constant. As the proof clarifies, this additional
dependence arises since the regularizer over the subspace $\overline
{\ModelSet}$ is
larger than the norm ${\|\cdot\|}$ by a factor of at most
${\Psi(\overline{\ModelSet})}$ (see Definition~\ref{DefnSubspaceCompat}).

Obtaining concrete rates using Corollary~\ref{CorExact} requires some
work in order to verify the conditions of Theorem~\ref{ThmMain} and
to provide control on the three quantities in the
bounds (\ref{EqnThreePart}) and (\ref{EqnThreePartReg}), as
illustrated in the examples to follow.

%%%%%%%%%%%%%%%%%%%%%%%%%%%%%%%%%%%%%%%%%%%%%%%%%%%%%%%%%%%%%%%%%%%%%%%%%%%%%%%%%%%%%%%%%%%%%%

%s4 #&#
\section{Convergence Rates for Sparse Regression}
\label{SecSparseReg}

As an illustration, we begin with one of the simplest statistical
models, namely, the standard linear model. It is based on $\numobs$
observations $Z_i = (x_i, y_i) \in\real^\pdim\times\real$ of
covariate-response pairs. Let $y \in\real^\numobs$ denote a vector
of the responses, and let $\Xmat\in\real^{\numobs\times\pdim}$ be
the design matrix, where $x_i \in\real^\pdim$ is the $i$th row.
This pair is linked via the linear model
%
%e26 #&#
\begin{equation}
\label{EqnLinRegGen} y = \Xmat{\betapar^*}+ w,
\end{equation}
where ${\betapar^*}\in\real^\pdim$ is the unknown regression vector
and $w \in\real^\numobs$ is a noise vector. To begin, we focus on
this simple linear setup and describe extensions to generalized
models in Section~\ref{SecGenLinModel}.

Given the data set ${Z_1^{n}}= (y, X) \in\real^\numobs\times
\real^{\numobs\times\pdim}$, our goal is to obtain a ``good''
estimate $\widehat{\betapar}$ of the regression vector ${\betapar^*}$, assessed
either in terms of its $\ell_2$-error $\|\widehat{\betapar}-
{\betapar^*}\|_2$ or
its $\ell_1$-error $\|\widehat{\betapar}- {\betapar^*}\|_1$. It is
worth noting
that whenever $\pdim> \numobs$, the standard linear
model (\ref{EqnLinRegGen}) is unidentifiable in a certain sense, since
the rectangular matrix $\Xmat\in\real^{\numobs\times\pdim}$ has a
null space of dimension at least $\pdim- \numobs$. Consequently, in
order to obtain an identifiable model---or at the very least, to bound
the degree of nonidentifiability---it is essential to impose
additional constraints on the regression vector ${\betapar^*}$. One
natural constraint is some type of sparsity in the regression vector;
for instance, one might assume that ${\betapar^*}$ has at most $\kdim$
nonzero coefficients, as discussed at more length in
Section~\ref{SecLassoHard}. More generally, one might assume that
although ${\betapar^*}$ is not exactly sparse, it can be
well-approximated by a sparse vector, in which case one might say that
${\betapar^*}$ is ``weakly sparse,'' ``sparsifiable'' or
``compressible.'' Section~\ref{SecWeakSparsity} is devoted to a more
detailed discussion of this weakly sparse case.

A natural $M$-estimator for this problem is the
Lasso \cite{Chen98,Tibshirani96}, obtained by solving the
$\ell_1$-penalized\vadjust{\goodbreak} quadratic program
%
%e27 #&#
\begin{equation}
\label{EqnLassoHardGen} {\widehat{\betapar}_{\regpar}}\in\arg
\min_{\param\in\real
^\pdim} \biggl\{ \frac{1}{2 \numobs} \|y - X \param\|_2^2
+ \regpar\|\param\|_1 \biggr\}
\end{equation}
for some choice $\regpar> 0$ of regularization parameter. Note that
this Lasso estimator is a particular case of the general
$M$-estimator (\ref{EqnMest}), based on the loss function and
regularization pair $\Loss(\betapar; {Z_1^{n}}) =
\frac{1}{2\numobs} \|y - X \param\|_2^2$ and ${\Regplain(\param)} =
\sum_{j=1}^\pdim|\param_j| = \|\param\|_1$. We now show how
Theorem~\ref{ThmMain} can be specialized to obtain bounds on the error
${\widehat{\betapar}_{\regpar}}- {\betapar^*}$ for the Lasso estimate.

%%%%%%%%%%%%%%%%%%%%%%%%%%%%%%%%%%%%%%%%%%%%%%%%%%%%%%%%%%%%%%%%%%%%%%%%%%%%%%%%%%%%%%%%%%%%%%%

%s4.1 #&#
\subsection{Restricted Eigenvalues for Sparse Linear Regression}
\label{SecRE}

For the least squares loss function that underlies the Lasso, the
first-order Taylor series expansion from Definition~\ref{DefnRSC} is
exact, so that
\[
\delLike\bigl(\delpar, \theta^{*}\bigr) = \biggl\langle{\delpar}, {
\frac
{1}{\numobs} \Xmat^T \Xmat\delpar} \biggr\rangle=
\frac
{1}{\numobs} \|\Xmat\delpar\|_2^2.
\]
Thus, in this special case, the Taylor series error is independent of
$\param^*$, a fact which allows for substantial theoretical
simplification. More precisely, in order to establish restricted
strong convexity, it suffices to establish a lower bound on $\|\Xmat
\delpar\|_2^2/\numobs$ that holds uniformly for an appropriately
restricted subset of $\pdim$-dimensional vectors $\delpar$.

As previously discussed in Example~\ref{ExaSparseVec}, for any subset
$\Sset\subseteq\{1, 2, \ldots, \pdim\}$, the $\ell_1$-norm is
decomposable with respect to the subspace $\ModelSet(\Sset) =
\{\TmpVarA\in\real^\pdim\mid\TmpVarA_{{{\Sset^c}}} = 0 \}$
and its
orthogonal complement. When the unknown regression vector $\param^*
\in\real^\pdim$ is exactly sparse, it is natural to choose $\Sset$
equal to the support set of ${\betapar^*}$. By appropriately
specializing the definition (\ref{EqnDefnConeSet}) of $\ConvSet$, we
are led to consider the cone
%
%e28 #&#
\begin{equation}
\label{EqnHardSparseCone} \ConvSet(\Sset) :=\bigl\{ \Delta\in
\real^\pdim\mid\|\Delta_{{{\Sset^c}}}\|_1 \leq3 \|
\Delta_\Sset\|_1 \bigr\}.
\end{equation}
See Figure~\ref{FigConeSet}(a) for an illustration of this set in
three dimensions. With this choice, restricted strong convexity with
respect to the $\ell_2$-norm is equivalent to requiring that the
design matrix $\Xmat$ satisfy the condition
%
%e29 #&#
\begin{equation}
\label{EqnRE} \frac{\|\Xmat\param\|_2^2}{\numobs} \geq\NewMcon\|\param
\|_2^2\quad
\mbox{for all $\param\in\ConvSet(\Sset)$.}
\end{equation}
This lower bound is a type of \textit{restricted eigenvalue} (RE)
condition and has been studied in past work on basis pursuit and the
Lasso (e.g., \cite{BiRiTsy08,MeiYu09,RasWaiYu09,GeerBuhl09}). One
could also require that a related condition hold with respect to the
$\ell_1$-norm---viz.
%
%e30 #&#
\begin{equation}
\label{EqnREOne} \frac{\|\Xmat\param\|_2^2}{\numobs} \geq\NewMcon'
\frac{\|\param\|_1^2}{|\Sset|} \quad\mbox{for all $\param\in\ConvSet(\Sset)$.}
\end{equation}
This type of $\ell_1$-based RE condition is less restrictive than the
corresponding $\ell_2$-version (\ref{EqnRE}). We refer the reader to
the paper by van de Geer and B\"{u}hlmann \cite{GeerBuhl09} for an
extensive discussion of different types of restricted eigenvalue or
compatibility conditions.

It is natural to ask whether there are many matrices that satisfy
these types of RE conditions. If $X$ has i.i.d. entries following a
sub-Gaussian distribution (including Gaussian and Bernoulli variables
as special cases), then known results in random matrix theory imply
that the restricted isometry property \cite{CandesTao06} holds with
high probability, which in turn implies that the RE condition
holds \cite{BiRiTsy08,GeerBuhl09}. Since statistical applications
involve design matrices with substantial dependency, it is natural to
ask whether an RE condition also holds for more general random
designs. This question was addressed by Raskutti et
al. \cite{RasWaiYu09,RasWaiYu10}, who showed that if the design matrix
$\Xmat\in\real^{\numobs\times\pdim}$ is formed by independently
sampling each row $X_i \sim N(0, \Sigma)$, referred to as the
\textit{$\Sigma$-Gaussian ensemble}, then there are strictly
positive constants $(\kappa_1, \GARTWO)$, depending only on the
positive definite matrix $\Sigma$, such that
%
%e31 #&#
\begin{eqnarray}
\label{EqnGarvesh}
\frac{\|\Xmat\param\|^2_2}{\numobs} &\geq&\kappa_1 \|
\param
\|_2^2 - \GARTWO\frac{\log\pdim}{\numobs} \|\param
\|_1^2 \nonumber\\[-8pt]\\[-8pt]
&&\eqntext{\mbox{for all $\param\in\real^\pdim$}}
\end{eqnarray}
with probability greater than $1 - c_1 \exp(-c_2 \numobs)$. The
bound (\ref{EqnGarvesh}) has an important consequence: it guarantees
that the RE property (\ref{EqnRE}) holds\footnote{To see this fact,
note that for any $\param\in\ConvSet(\Sset)$, we have
$\|\param\|_1 \leq4 \|\param_\Sset\|_1 \leq4 \sqrt{\kdim}
\|\param_\Sset\|_2$. Given the lower bound (\ref{EqnGarvesh}), for
any $\param\in\ConvSet(\Sset)$, we have\vspace*{-2pt} the lower bound
$\frac{\|\Xmat\param\|_2}{\sqrt{\numobs}} \geq\{ \kappa_1
- 4 \kappa_2 \sqrt{\frac{\kdim\log\pdim}{\numobs}} \}
\|\param\|_2 \geq\frac{\kappa_1}{2} \|\param\|_2$, where
final inequality follows as long as $\numobs> 64
(\kappa_2/\kappa_1)^2 \kdim\log\pdim$.} with $\NewMcon=
\frac{\kappa_1}{2} > 0$ as long as $\numobs> 64 (\kappa_2/\kappa_1)
\kdim\log\pdim$. Therefore, not only do there exist matrices
satisfying the RE property (\ref{EqnRE}), but any matrix sampled from
a $\Sigma$-Gaussian ensemble will satisfy it with high
probability. Related analysis by Rudelson and Zhou \cite{RudZho11}
extends these types of guarantees to the case of sub-Gaussian designs,
also allowing for substantial dependencies among the covariates.

%%%%%%%%%%%%%%%%%%%%%%%%%%%%%%%%%%%%%%%%%%%%%%%%%%%%%%%%%%%%%%%%%%%%%%%%%%%%%%%%%%%%%%%%%%%%%%

%s4.2 #&#
\subsection{Lasso Estimates with Exact Sparsity}
\label{SecLassoHard}

We now show how Corollary~\ref{CorExact} can be used to derive
convergence rates for the error of the Lasso estimate when the unknown
regression vector $\param^*$ is $\kdim$-sparse. In order to state
these results, we require some additional notation.\vadjust{\goodbreak} Using $X_j
\in\real^\numobs$ to denote the $j$th column of~$\Xmat$, we say
that $\Xmat$ is \textit{column-normalized} if
%
%e32 #&#
\begin{equation}
\label{EqnColNorm} \frac{\|X_j\|_2}{\sqrt{\numobs}} \leq1 \quad\mbox{for
all $j = 1, 2, \ldots,
\pdim$.}
\end{equation}
Here we have set the upper bound to one in order to simplify notation.
This particular choice entails no loss of generality, since we can
always rescale the linear model appropriately (including the
observation noise variance) so that it holds.

In addition, we assume that the noise vector $w \in\real^\numobs$ is
zero-mean and has \textit{sub-Gaussian tails}, meaning that there is a
constant $\sigma> 0$ such that for any fixed \mbox{$\|v\|_2 = 1$},
%
%e33 #&#
\begin{equation}
\label{EqnSubGaussNoise} \mprob\bigl[ \bigl|\langle{v}, {w} \rangle\bigr| \geq t
\bigr]
\leq2 \exp\biggl( - \frac{\delta^2}{2 \sigma^2} \biggr) \quad\mbox{for all
$\delta> 0$.}\hspace*{-27pt}
\end{equation}
For instance, this condition holds when the noise vector $w$ has
i.i.d. $N(0,1)$ entries or consists of independent bounded random
variables. Under these conditions, we recover as a corollary of
Theorem~\ref{ThmMain} the following result:
%
%co2 #&#
\begin{cors}
\label{CorLassoHard}
Consider an $\kdim$-sparse instance of the linear regression
model (\ref{EqnLinRegGen}) such that $\Xmat$ satisfies the RE
condition (\ref{EqnRE}) and the column normalization
condition (\ref{EqnColNorm}). Given the Lasso
program (\ref{EqnLassoHardGen}) with regularization parameter
$\regpar= 4 \sigma\sqrt{\frac{ \log\pdim}{\numobs}}$,
then with probability at least $1 - \plaincon_1 \exp(-\plaincon_2
\numobs\regpar^2)$, any\vspace*{1pt} optimal solution ${\widehat{\betapar
}_{\regpar}}$ satisfies the
bounds
%
%e34 #&#
\begin{eqnarray}
\label{EqnLassoErrHard}
\bigl\|{\widehat{\betapar}_{\regpar}}- {\betapar^*}
\bigr\|^2_2 &\leq&\frac
{64
\sigma^2}{\NewMcon^2} \frac{\kdim\log\pdim}{\numobs}\quad
\mbox{and}\nonumber\\[-8pt]\\[-8pt]
\bigl\|{\widehat{\betapar}_{\regpar}}- {\betapar^*}
\bigr\|_1 &\leq&\frac
{24 \sigma}{\NewMcon} \kdim\sqrt{\frac{\log\pdim}{\numobs}}.\nonumber
\end{eqnarray}
\end{cors}

Although error bounds of this form are known from past work
(e.g., \cite{BiRiTsy08,CandesTao06,MeiYu09}), our proof illuminates
the underlying structure that leads to the different terms in the
bound---in particular, see equations (\ref{EqnThreePart})
and (\ref{EqnThreePartReg}) in the statement of
Corollary~\ref{CorExact}.
\begin{pf*}{Proof of Corollary~\ref{CorLassoHard}}
We first note that the RE condition (\ref{EqnREOne}) implies that RSC
holds with respect to the subspace $\ModelSet(\Sset)$. As discussed
in\break \mbox{Example}~\ref{ExaSparseVec}, the $\ell_1$-norm is
decomposable with respect to $\ModelSet(\Sset)$ and its orthogonal
complement, so that we may set $\overline{\ModelSet}(\Sset) =
\ModelSet(\Sset)$.
Since any vector $\theta\in\ModelSet(\Sset)$ has at most $\kdim$
nonzero entries, the subspace compatibility constant is given by
${\Psi(\ModelSet(\Sset))} =\break \sup_{\theta\in
\ModelSet(\Sset) \setminus\{0\}} \frac{\|\theta\|_1}{\|\theta\|_2}
= \sqrt{\kdim}$.

The final step is to compute an appropriate choice of the
regularization parameter. The gradient of the quadratic loss is given
by $\nabla\Loss(\betapar; (y,X)) = \frac{1}{\numobs} \Xmat^T
w$, whereas the dual norm of the $\ell_1$-norm is the
$\ell_\infty$-norm. Consequently, we need to specify a choice of
$\regpar> 0$ such that
\[
\regpar\geq2 {\Regplain^*\bigl(\nabla\Loss\bigl({\betapar^*}\bigr)\bigr
)} = 2 \biggl\|
\frac{1}{\numobs} \Xmat^T w \biggr\|_\infty
\]
with high probability. Using the column
normalization (\ref{EqnColNorm}) and
sub-Gaussian (\ref{EqnSubGaussNoise}) conditions, for each $j = 1,
\ldots, \pdim$, we have the tail bound $\mprob[
|\langle{X_j}, {w} \rangle/\numobs| \geq t ] \leq2 \exp( -
\frac{\numobs t^2}{2 \sigma^2} )$. Consequently, by union bound,
we conclude that $\mprob[ \|\Xmat^T w/\numobs\|_\infty\geq t
] \leq2 \exp( - \frac{\numobs t^2}{2 \sigma^2} +\break \log
\pdim)$. Setting $t^2 = \frac{4 \sigma^2 \log
\pdim}{\numobs}$, we see that the choice of $\regpar$ given in
the statement is valid with probability at least $1 - c_1 \exp(-c_2
\numobs\regpar^2)$. Consequently, the claims (\ref{EqnLassoErrHard})
follow from the bounds (\ref{EqnThreePart})
and (\ref{EqnThreePartReg}) in Corollary~\ref{CorExact}.
\end{pf*}
%

%%%%%%%%%%%%%%%%%%%%%%%%%%%%%%%%%%%%%%%%%%%%%%%%%%%%%%%%%%%%%%%%%%%%%%%%%%%%%%%%%%%%%%%%%%%%%%%

%s4.3 #&#
\subsection{Lasso Estimates with Weakly Sparse Models}
\label{SecWeakSparsity}

We now consider regression models for which $\param^*$ is not
exactly sparse, but rather can be approximated well by a sparse
vector. One way in which to formalize this notion is by considering
the $\ell_\qpar$ ``ball'' of radius $\radq$, given by
\begin{eqnarray}
\Ball_\qpar(\radq) &:=&\Biggl\{ \param\in\real^\pdim\Bigm|
\sum_{i=1}^\pdim|\param_i|^\qpar
\leq\radq\Biggr\}\nonumber\\
&&\eqntext{\mbox{where $\qpar\in[0,1]$ is fixed.}}
\end{eqnarray}
In the special case $\qpar= 0$, this set corresponds to an exact
sparsity constraint---that is, $\theta^{*}\in\Ball_0(R_0)$ if
and only if $\theta^{*}$ has at most $R_0$ nonzero entries. More
generally, for $\qpar\in(0,1]$, the set $\Ball_\qpar(\radq)$
enforces a certain decay rate on the ordered absolute values of~$\param^*$.

In the case of weakly sparse vectors, the constraint set $\ConeSet$
takes the form
%
%e35 #&#
\begin{eqnarray}
\label{EqnConeSetWeak}
\hspace*{17pt}&&\ConeSet\bigl(\ModelSet, \overline{\ModelSet};
\param^*
\bigr) \nonumber\\[-8pt]\\[-8pt]
&&\quad= \bigl\{ \delpar\in\real^\pdim\mid\|\delpar_{{\Sset^c}}
\|_1 \leq3 \|\delpar_\Sset\|_1 + 4 \bigl\|
\param^*_{{{\Sset^c}}}\bigr\|_1 \bigr\}.\nonumber
\end{eqnarray}
In contrast to the case of exact sparsity, the set $\ConeSet$ is no
longer a cone, but rather contains a ball centered at the
origin---compare panels (a) and (b) of Figure~\ref{FigConeSet}. As a
consequence, it is \textit{never} possible to ensure that $\|\Xmat
\theta\|_2/\sqrt{\numobs}$ is uniformly bounded from below for all
vectors $\theta$ in the set (\ref{EqnConeSetWeak}), and so a strictly
positive tolerance term $\rscslop(\theta^{*}) > 0$ is required. The
random matrix result (\ref{EqnGarvesh}), stated in the previous
section, allows us to establish a form of RSC that is appropriate for
the setting of $\ell_\qpar$-ball sparsity. We summarize our conclusions
in the following:
%
%co3 #&#
\begin{cors}
\label{CorLassoWeak}
Suppose that $\Xmat$ satisfies the RE condition (\ref{EqnGarvesh}) as
well as the column normalization condition (\ref{EqnColNorm}), the
noise $w$ is sub-Gaussian (\ref{EqnSubGaussNoise}) and $\param^*$
belongs to $\Ball_\qpar(\radq)$ for a radius $\radq$ such that
$\sqrt{\radq} (\frac{\log\pdim}{\numobs} )^{
{1}/{2} -
{\qpar}/{4}} \leq1$. Then if we solve the Lasso with
regularization parameter $\regpar= 4 \sigma\sqrt{\frac{\log
\pdim}{\numobs}}$, there are universal positive constants
$(\plaincon_0, \plaincon_1, \plaincon_2)$ such that any optimal
solution ${\widehat{\betapar}_{\regpar}}$ satisfies
%
%e36 #&#
\begin{equation}
\label{EqnWeakUpperBound} \bigl\|{\widehat{\betapar}_{\regpar}}- {\betapar^*}
\bigr\|^2_2 \leq c_0 \radq\biggl(
\frac{\sigma^2}{\kappa_1^2} \frac{\log\pdim}{\numobs} \biggr)^{1-
{\qpar}/{2}}
\end{equation}
with probability at least $1 - \plaincon_1 \exp(-\plaincon_2 \numobs
\regpar^2)$.
\end{cors}

%pa4.3.subsubsection.1 #&#
\begin{Remarks*} Note that this corollary is a strict
generalization of Corollary~\ref{CorLassoHard}, to which it reduces
when $\qpar= 0$. More generally, the parameter $\qpar\in[0,1]$
controls the relative ``sparsifiability'' of $\param^*$, with lar\-ger
values corresponding to lesser sparsity. Naturally then, the rate
slows down as $\qpar$ increases from $0$ toward $1$. In fact,
Raskutti et al. \cite{RasWaiYu09} show that the
rates (\ref{EqnWeakUpperBound}) are minimax-optimal over the
$\ell_\qpar$-balls---implying that not only are the consequences of
Theorem~\ref{ThmMain} sharp for the Lasso, but, more generally, no
algorithm can achieve faster rates.
\end{Remarks*}
\begin{pf*}{Proof of Corollary~\ref{CorLassoWeak}}
Since the loss function $\Loss$ is quadratic, the proof of Corollary
\ref{CorLassoHard} shows that the stated choice $\regpar= 4
\sqrt{\frac{ \sigma^2 \log\pdim}{\numobs}}$ is valid\vspace*{1pt} with
probability at least $1-c \exp(-c' \numobs\regpar^2)$. Let us now show
that the RSC condition holds. We do so via condition (\ref{EqnGarvesh})
applied to equation (\ref{EqnConeSetWeak}). For a threshold $\taupar>
0$ to be chosen, define the thresholded subset
%
%e37 #&#
\begin{equation}
\label{EqnDefnStau} \Stau:=\bigl\{ j \in\{1, 2, \ldots, \pdim\}
\mid\bigl|
\betapar^*_j\bigr| > \taupar\bigr\}.
\end{equation}
Now recall the subspaces $\ModelSet(\Stau)$ and ${\ModelSet^\perp
}(\Stau)$
previously defined in equations (\ref{EqnModelSparse})
and (\ref{EqnBsetSparse}) of Example~\ref{ExaSparseVec}, where we set
$\Sset= \Stau$. The following lemma, proved in the
supplement \cite{Neg12supp}, provides sufficient conditions for
restricted strong convexity with respect to these subspace pairs:
%
%le2 #&#
\begin{lems}
\label{LemRSCWeak}
Suppose that the conditions of Corollary~\ref{CorLassoWeak} hold and
$\numobs> 9 \kappa_2 |\Stau| \log\pdim$. Then with the
choice
$\taupar
= \frac{\regpar}{\kappa_1}$, the RSC condition holds over
$\ConvSet(\ModelSet(\Stau),\break {\ModelSet^\perp}(\Stau),
\param^*)$ with $\NewMcon= \kappa_1/4$ and $\rscslop^2 = 8
\kappa_2
\frac{\log\pdim}{\numobs} \|\theta^{*}_{\Sset_\taupar^c} \|_1^2$.
\end{lems}

Consequently, we may apply Theorem~\ref{ThmMain} with
$\NewMcon= \kappa_1/4$ and $\rscslop^2(\theta^{*}) = 8 \kappa_2
\frac{\log\pdim}{\numobs} \|\theta^{*}_{\Sset_\taupar^c}\|_1^2$
to conclude
that
%
%e38 #&#
\begin{eqnarray}
\label{EqnEvans}
&&\bigl\| {\widehat{\betapar}_{\regpar}}- {\betapar^*}
\bigr\|_2^2 \nonumber\\
&&\quad\leq144 \frac{\regpar^2}{\kappa_1^2} |\Stau| \\
&&\qquad{}+
\frac{4
\regpar}{\kappa_1} \biggl\{ 16 \kappa_2 \frac{\log\pdim}{\numobs} \bigl\|
\theta^{*}_{\Sset_\taupar^c}\bigr\|_1^2 + 4 \bigl\|
\theta^{*}_{\Sset_\taupar^c}\bigr\|_1 \biggr\},\nonumber
\end{eqnarray}
where we have used the fact that ${\Psi^2(\Stau)} = |\Stau|$, as
noted in the proof of Corollary~\ref{CorLassoHard}.

It remains to upper bound the cardinality of $\Stau$ in terms of the
threshold $\taupar$ and $\ell_\qpar$-ball radius $\radq$. Note that we
have
%
%e39 #&#
\begin{equation}
\label{EqnUpperStau} \radq\geq\sum_{j =1}^\pdim\bigl|
\theta^{*}_j\bigr|^q \geq\sum
_{j \in
\Stau} \bigl|\theta^{*}_i\bigr|^q
\geq\taupar^\qpar|\Stau|,
\end{equation}
hence, $|\Stau| \leq\taupar^{-\qpar} \radq$ for any $\taupar>
0$. Next
we upper bound the approximation error $\|\theta^{*}_{\Sset_\taupar^c}\|_1$,
using the fact that $\theta^{*}\in\Ball_\qpar(\radq)$. Letting
$\Sset_\taupar^c$ denote the complementary set $\Stau\setminus\{1, 2,
\ldots, \pdim\}$, we have
%
%e40 #&#
\begin{eqnarray}
\label{EqnUpperApp} \bigl\|\theta^{*}_{\Sset_\taupar^c}\bigr\|_1 &=&
\sum_{j \in{\Sset_\taupar^c}} \bigl|\theta^{*}_j\bigr| =
\sum_{j \in{\Sset_\taupar^c}} \bigl|\theta^{*}_j\bigr|^\qpar\bigl|
\theta^{*}_j\bigr|^{1-\qpar} \nonumber\\[-8pt]\\[-8pt]
&\leq&\radq
\taupar^{1-\qpar}.\nonumber
\end{eqnarray}
Setting $\taupar= \regpar/\kappa_1$ and then substituting the
bounds (\ref{EqnUpperStau}) and (\ref{EqnUpperApp}) into the
bound (\ref{EqnEvans}) yields
\begin{eqnarray*}
\bigl\| {\widehat{\betapar}_{\regpar}}- {\betapar^*}\bigr\|_2^2
% \radq+ 4 \big(\frac{\regpar^2}{\kappa_1^2} \big)^{1 -
% \frac{\qpar}{2}} \radq\bigg\{ 16 \kappa_2 \frac{\log
% \pdim}{\numobs} \radq\big(\frac{\regpar}{\kappa_1}
% \big)^{1-\qpar} + 4 \bigg\} \\
&\leq&160 \biggl(
\frac{\regpar^2}{\kappa_1^2} \biggr)^{1 -
{\qpar}/{2}} \radq\\
&&{}+ 64 \kappa_2 \biggl\{
\biggl(\frac{\regpar^2}{\kappa_1^2} \biggr)^{1 - {\qpar}/{2} }
\radq\biggr\}^2
\frac{(\log\pdim)/\numobs}{\regpar/\kappa_1}.
\end{eqnarray*}
For any fixed noise variance, our choice of regularization parameter
ensures that the ratio $\frac{(\log\pdim)/\numobs}{\regpar/\kappa_1}$
is of order one, so that the claim follows.
\end{pf*}
%

%%%%%%%%%%%%%%%%%%%%%%%%%%%%%%%%%%%%%%%%%%%%%%%%%%%%%%%%%%%%%%%%%%%%%%%%%%%%%%%%%%%%%%%%%%%%%%%

%s4.4 #&#
\subsection{Extensions to Generalized Linear Models}
\label{SecGenLinModel}

In this section we briefly outline extensions of the preceding
results to the family of generalized linear models (GLM). Suppose
that conditioned on a vector $x \in\real^\pdim$ of covariates, a
response variable $y \in\mathcal{Y}$ has the distribution
%
%e41 #&#
\begin{equation}
\label{EqnGLM}\quad \mprob_{\param^*}(y \mid x) \propto\exp\biggl\{
\frac{ y
\langle{\param^*}, {x} \rangle- \Phi(\langle{\param^*}, {x}
\rangle)}{c(\sigma)} \biggr\}.
\end{equation}
Here the quantity $c(\sigma)$ is a fixed and known scale parameter,
and the function $\Phi\dvtx\real\rightarrow\real$ is the link function,
also known. The family (\ref{EqnGLM}) includes many well-known
classes of regression models as special cases, including ordinary
linear regression [obtained with $\mathcal{Y} = \real$, $\Phi(t) =
t^2/2$ and $c(\sigma) = \sigma^2$] and logistic regression [obtained
with $\mathcal{Y} = \{0,1\}$, $c(\sigma) = 1$ and $\Phi(t) = \log(1 +
\exp(t))$].

Given samples $Z_i = (x_i, y_i) \in\real^\pdim\times\mathcal{Y}$,
the goal is to estimate the unknown vector $\param^*\in
\real^\pdim$. Under a sparsity assumption on $\param^*$, a natural
estimator is based on minimizing the (negative) log likelihood,
combined with an $\ell_1$-regularization term. This combination leads
to the convex program
%
%e42 #&#
\begin{eqnarray}
\label{EqnGLMEst}\hspace*{23pt} {\widehat{\theta}}_{\regpar} &\in&\arg\min_{\theta\in
\real^\pdim
}
\Biggl\{ \underbrace{\frac{1}{\numobs} \sum_{i=1}^\numobs
\bigl\{ -y_i \langle{\theta}, {x_i} \rangle+ \Phi\bigl(
\langle{\theta}, {x_i} \rangle\bigr) \bigr\}}_{\Loss(\theta;
{Z_1^{n}})} \nonumber\\[-8pt]\\[-8pt]
&&\hspace*{135.2pt}{}+
\regpar\|\theta\|_1 \Biggr\}.\nonumber
\end{eqnarray}
In order to extend the error bounds from the previous section, a key
ingredient is to establish that this GLM-based loss function satisfies
a form of restricted strong convexity. Along these lines, Negahban et
al. \cite{Neg09} proved the following result: suppose that the
covariate vectors $x_i$ are zero-mean with covariance matrix $\Sigma
\succ0$ and are drawn i.i.d. from a distribution with sub-Gaussian
tails [see equation (\ref{EqnSubGaussNoise})]. Then there are
constants $\kappa_1, \kappa_2$ such that the first-order Taylor series
error for the GLM-based loss (\ref{EqnGLMEst}) satisfies the lower
bound
%
%e43 #&#
\begin{eqnarray}
\label{EqnSahandGLM} \delLike\bigl(\delpar, \theta^{*}\bigr) &\geq&
\kappa_1 \|\delpar\|_2^2 -
\kappa_2 \frac{\log\pdim}{\numobs} \|\delpar\|_1^2\nonumber\\[-8pt]\\[-8pt]
&&\eqntext{\mbox{for all $\|\delpar\|_2 \leq1$.}}
\end{eqnarray}
As discussed following Definition~\ref{DefnRSC}, this type of lower
bound implies that $\Loss$ satisfies a form of RSC, as long as the
sample size scales as $\numobs= \Omega(\kdim\log\pdim)$, where
$\kdim$ is the target sparsity. Consequently, this lower
bound (\ref{EqnSahandGLM}) allows us to recover analogous bounds on
the error $\|{\widehat{\theta}}_{\regpar} - \param^*\|_2$ of the GLM-based
estimator (\ref{EqnGLMEst}).

%s5 #&#
\section{Convergence Rates for Group-Structured Norms}
\label{SecGroup}

The preceding two sections addressed $M$-estima\-tors based on
$\ell_1$-regularization, the simplest type of decomposable
regularizer. We now turn to some extensions of our results to more
complex regularizers that are also decomposable. Various researchers
have proposed extensions of the Lasso based on regularizers that have
more structure than the $\ell_1$-norm
(e.g., \cite{Turlach05,YuaLi06,ZhaRoc06,Mei08,Bar08}). Such
regularizers allow one to impose different types of block-sparsity
constraints, in which groups of parameters are assumed to be active
(or inactive) simultaneously. These norms arise in the context of
\mbox{multivariate} regression, where the goal is to predict a multivariate
output in $\real^\gsize$ on the basis of a set of $\pdim$ covariates.
Here it is appropriate to assume that groups of covariates are useful
for predicting the different elements of the $\gsize$-dimensional
output vector. We refer the reader to the
papers \cite{Turlach05,YuaLi06,ZhaRoc06,Mei08,Bar08} for further
discussion of and motivation for the use of block-structured norms.

Given a collection $\GROUP= \{{G}_1, \ldots, {G}_\numgroup\}$
of groups, recall from Example~\ref{ExaBlockSparse} in
Section~\ref{DecompSection} the definition of the group norm
${\|\cdot\|_{\GROUP, \SUPERGVEC}}$. In full generality, this group norm
is based on a weight vector $\SUPERGVEC= (\gqpar_1, \ldots,
\gqpar_\numgroup) \in[2, \infty]^{\numgroup}$, one for each group.
For simplicity, here we consider the case when $\gqpar_t =
\gqpar$ for all $t = 1, 2, \ldots, \numgroup$, and we use
${\|\cdot\|_{\GROUP, \gqpar}}$ to denote the associated group norm.
As a
natural extension of the Lasso, we consider the \textit{block Lasso}
estimator
%
%e44 #&#
\begin{equation}\hspace*{5pt}
{\widehat{\theta}}\in\arg\min_{\betapar\in\real^\pdim} \biggl\{ \frac
{1}{\numobs} \|y - X
\betapar\|_2^2 + \regpar{\|\betapar\|_{\GROUP, \SUPERGVEC}}
\biggr\},
\end{equation}
where $\regpar> 0$ is a user-defined regularization parameter.
Different choices of the parameter $\gqpar$ yield different
estimators, and in this section we consider the range $\gqpar\in[2,
\infty]$. This range covers the two most commonly applied choices,
$\gqpar= 2$, often referred to as the group Lasso, as well as the
choice $\gqpar= +\infty$.

%s5.1 #&#
\subsection{Restricted Strong Convexity for Group Sparsity}

As a parallel to our analysis of ordinary sparse regression, our first
step is to provide a condition sufficient to guarantee restricted
strong convexity for the group-sparse setting. More specifically, we
state the natural extension of condition (\ref{EqnGarvesh}) to the
block-sparse setting and prove that it holds with high probability
for the class of $\Sigma$-Gaussian random designs. Recall from
Theorem~\ref{ThmMain} that the dual norm of the regularizer plays a
central role. As discussed previously, for the block-$(1,
\gqvec)$-regularizer, the associated dual norm is a block-$(\infty,
{{\gqpar^\ast}})$ norm, where $(\gqpar, {{\gqpar^\ast}})$ are conjugate
exponents satisfying \mbox{$\frac{1}{\gqpar} + \frac{1}{{{\gqpar^\ast}}}
= 1$}.

Letting $\basgauss\sim N(0, I_{\pdim\times\pdim})$ be a standard
normal vector, we consider the following condition. Suppose that
there are strictly positive constants\vadjust{\goodbreak} $(\kappa_1, \kappa_2)$ such
that, for all $\delpar\in\real^\pdim$, we have
%
%e45 #&#
\begin{equation}
\label{EqnREGroup} \frac{\|\Xmat\delpar\|^2_2}{\numobs} \geq\kappa_1 \|
\delpar
\|^2_2 - \kappa_2 \DualNoise^2
\bigl({{\gqpar^\ast}}\bigr) \|\delpar\|^2_{1,
\gqparbar},
\end{equation}
where $\DualNoise({{\gqpar^\ast}}) :=\Exs[ \max_{t=1,2,
\ldots,
\numgroup}\frac{\|\basgauss_{{G}_t}\|_{{{\gqpar^\ast}}}}{\sqrt{\numobs
}} ]$.
To understand this condition, first consider the special case of
$\numgroup= \pdim$ groups, each of size one, so that the group-sparse
norm reduces to the ordinary $\ell_1$-norm, and its dual is the
$\ell_\infty$-norm. Using $\gqpar= 2$ for concreteness, we have
$\DualNoise(2) = \Exs[ \|\basgauss\|_\infty]/\sqrt{\numobs}
\leq\sqrt{ \frac{3 \log\pdim}{\numobs}}$ for all $\pdim\geq10$,
using standard bounds on Gaussian maxima. Therefore,
condition (\ref{EqnREGroup}) reduces to the earlier
condition (\ref{EqnGarvesh}) in this special case.

Let us consider a more general setting, say, with \mbox{$\gqpar= 2$} and
$\numgroup$ groups each of size $\Gmax$, so that $\pdim= \numgroup
\Gmax$. For this choice of groups and norm, we have
\[
\DualNoise(2) = \Exs\biggl[\max_{t = 1, \ldots, \numgroup} \frac{\|
\basgauss_{{G}_t}\|_2}{\sqrt{\numobs}} \biggr],
\]
where each sub-vector $w_{{G}_t}$ is a standard Gaussian vector
with $\Gmax$ elements. Since $\Exs[\|\basgauss_{{G}_t}\|_2] \leq
\sqrt{\Gmax}$, tail bounds for $\chi^2$-variates yield $\DualNoise(2)
\leq\sqrt{\frac{\Gmax}{\numobs}} + \sqrt{\frac{3 \log
\numgroup}{\numobs}}$, so that the condition (\ref{EqnREGroup}) is
equivalent to
%
%e46 #&#
\begin{eqnarray}
\label{EqnREGroupGeneral} \frac{\|\Xmat\delpar\|^2_2}{\numobs} &\geq&
\kappa_1 \|\delpar
\|^2_2 \nonumber\\
&&{}- \kappa_2 \biggl[\sqrt{
\frac{\Gmax}{\numobs}} + \sqrt{ \frac{3
\log
\numgroup}{\numobs}} \biggr]^2 {\|\delpar
\|_{\GROUP, 2}}^2 \\
&&\eqntext{\mbox{for all $\delpar\in\real^\pdim$.}}
\end{eqnarray}

Thus far, we have seen the form that condition (\ref{EqnREGroup})
takes for different choices of the groups and parameter~$\gqvec$. It
is natural to ask whether there are any matrices that satisfy the
condition (\ref{EqnREGroup}). As shown in the following result, the
answer is affirmative---more strongly, \mbox{almost} every matrix satisfied
from the $\Sigma$-Gaussian ensemble will satisfy this condition with
high probability. [Here we recall that for a nondegenerate
covariance matrix, a random design matrix $\Xmat\in\real^{\numobs
\times\pdim}$ is drawn from the $\Sigma$-Gaussian ensemble if each
row $x_i \sim N(0, \Sigma)$, i.i.d. for $i = 1, 2, \ldots, \numobs$.]
%
%pr1 #&#
\begin{props}
\label{PropREGroup}
For a design matrix $\Xmat\in\real^{\numobs\times\pdim}$ from the
$\Sigma$-ensemble, there are constants $(\kappa_1, \kappa_2)$
depending only on $\Sigma$ such that condition (\ref{EqnREGroup}) holds
with probability greater than $1 - c_1 \exp(-c_2 \numobs)$.
\end{props}

We provide the proof of this result in the
supplement~\cite{Neg12supp}. This condition can be used to show that
appropriate forms of RSC hold, for both the cases of exactly
group-sparse and weakly sparse vectors.\vadjust{\goodbreak} As with
$\ell_1$-regularization, these RSC conditions are milder than
analogous group-based RIP conditions
(e.g., \cite{HuaZha09,StohParHas09,Bar08}), which require that all
submatrices up to a certain size are close to isometries.

%s5.2 #&#
\subsection{Convergence Rates}

Apart from RSC, we impose one additional condition on the design
matrix. For a given group ${G}$ of size~$\Gmax$, let us view the
matrix $\Xmat_{{G}} \in\real^{\numobs\times\Gmax}$ as an
operator from $\ell^{\Gmax}_\gqpar\rightarrow\ell^\numobs_2$ and
define the associated operator norm $\norm\Xmat_{{G}} \norm_{{\gqpar
\rightarrow2}} :=\max_{\|\theta\|_\gqpar= 1}
\|\Xmat_{G} \theta\|_2$. We then require that
%
%e47 #&#
\begin{equation}\quad
\label{EqnColBlock} \frac{\norm\Xmat_{{G}_t} \norm_{{\gqpar\rightarrow
2}}}{\sqrt{\numobs}} \leq1 \quad\mbox{for all $t = 1, 2, \ldots,
\numgroup$.}
\end{equation}
Note that this is a natural generalization of the column normalization
condition (\ref{EqnColNorm}), to which it reduces when we have
$\numgroup= \pdim$ groups, each of size one. As before, we may
assume without loss of generality, rescaling $\Xmat$ and the noise as
necessary, that condition (\ref{EqnColBlock}) holds with constant one.
Finally, we define the maximum group size $\Gmax= \max_{t =
1, \ldots, \numgroup} |{G}_t|$. With this notation, we have the
following novel result:
%
%co4 #&#
\begin{cors}
\label{CorBlockNorm}
Suppose that the noise $w$ is sub-Gaussian (\ref{EqnSubGaussNoise}),
and the design matrix $\Xmat$ satisfies condition (\ref{EqnREGroup})
and the block normalization condition (\ref{EqnColBlock}). If we
solve the group Lasso with
%
%e48 #&#
\begin{equation}
\label{EqnBlockRegpar} \regpar\geq2 \sigma\biggl\{\frac{\Gmax^{1 -
1/\gqpar}}{\sqrt{\numobs}} + \sqrt{
\frac{\log\numgroup}{\numobs}} \biggr\},
\end{equation}
then with probability at least $1- 2/\numgroup^2$, for any group
subset $\GroupSset\subseteq\{1, 2, \ldots, \numgroup\}$ with
cardinality\break $|\GroupSset| = \groupcard$, any optimal solution
${\widehat{\betapar}_{\regpar}}$ satisfies
%
%e49 #&#
\begin{equation}
\label{EqnBlockqBound}\quad \bigl\|{\widehat{\betapar}_{\regpar}}-
\theta^{*}\bigr\|^2_2 \leq\frac{4
\regpar^2}{\NewMcon^2}
\groupcard+ \frac{4 \regpar}{\NewMcon} \sum_{t \notin
\GroupSset} \bigl\|
\param^*_{{G}_t} \bigr\|_\gqpar.
\end{equation}
\end{cors}
\begin{Remarks*} Since the result applies to any $\gqpar\in
[2, \infty]$, we can observe how the choices of different
group-sparse norms affect the convergence rates. So as to
simplify this discussion, let us assume that the groups are
all of equal size $\Gmax$, so that $\pdim= \Gmax\numgroup$
is the ambient dimension of the problem.

\textit{Case $\gqpar= 2$}: The case $\gqpar= 2$ corresponds
to the block $(1,2)$ norm, and the resulting estimator is frequently
referred to as the group Lasso. For this case, we can set the
regularization parameter\vspace*{-2pt} as $\regpar= 2 \sigma\{
\sqrt{\frac{\Gmax}{\numobs}} + \sqrt{\frac{\log\numgroup
}{\numobs}}
\}$. If we assume, moreover, that $\param^*$ is exactly
group-sparse, say, supported on a group subset $\GroupSset\subseteq
\{1, 2, \ldots, \numgroup\}$ of cardinality\vadjust{\goodbreak} $\groupcard$, then the
bound (\ref{EqnBlockqBound}) takes the form
%
%e50 #&#
\begin{equation}
\bigl\|\widehat{\betapar}- \theta^{*}\bigr\|_2^2
\precsim\frac
{\groupcard
\Gmax}{\numobs} + \frac{\groupcard\log\numgroup}{\numobs}.
\end{equation}
Similar bounds were derived in independent work by Lounici et
al. \cite{Lou09} and Huang and Zhang \cite{HuaZha09} for this special
case of exact block sparsity. The analysis here shows how the
different terms arise, in particular, via the noise magnitude measured
in the dual norm of the block regularizer.

In the more general setting of weak block sparsity,
Corollary~\ref{CorBlockNorm} yields a number of novel results. For
instance, for a given set of groups $\GROUP$, we can consider the
block sparse analog of the $\ell_\qpar$-``ball''---namely, the
set
\[
\Ball_\qpar(\radq; \GROUP, 2) :=\Biggl\{ \theta\in
\real^\pdim\Bigm|\sum_{t=1}^\numgroup\|
\theta_{{G}_t}\|_{2}^\qpar\leq\radq\Biggr\}.
\]
In this case, if we optimize the choice of $\Sset$ in the
bound (\ref{EqnBlockqBound}) so as to trade off the estimation and
approximation errors, then we obtain
\[
\bigl\|{\widehat{\theta}}- \param^*\bigr\|_2^2 \precsim\radq
\biggl( \frac{\Gmax}{\numobs} + \frac{\log\numgroup}{\numobs} \biggr
)^{1 -
{\qpar}/{2}},
\]
which is a novel result. This result is a generalization of our
earlier Corollary~\ref{CorLassoWeak}, to which it reduces when we have
$\numgroup= \pdim$ groups each of size $\Gmax= 1$.

\textit{Case $\gqpar= +\infty$}: Now consider the case of
$\ell_1/\ell_\infty$-regularization, as suggested in past
work \cite{Turlach05}. In this case, Corollary~\ref{CorBlockNorm}
implies that $\|\widehat{\betapar}- \theta^{*}\|^2_2 \precsim\frac
{\kdim
\Gmax^2 }{\numobs} + \frac{\kdim\log\numgroup}{\numobs}$. Similar
to the case $\gqpar= 2$, this bound consists of an estimation term
and a search term. The estimation term $\frac{\kdim
\Gmax^2}{\numobs}$ is larger by a factor of $\Gmax$, which
corresponds to the amount by which an $\ell_\infty$-ball in $\Gmax$
dimensions is larger than the corresponding $\ell_2$-ball.

We provide the proof of Corollary~\ref{CorBlockNorm} in the
supplementary appendix \cite{Neg12supp}. It is based on verifying the
conditions of Theorem~\ref{ThmMain}: more precisely, we use
Proposition~\ref{PropREGroup} in order to establish RSC, and we
provide a lemma that shows that the regularization
choice (\ref{EqnBlockRegpar}) is valid in the context of
Theorem~\ref{ThmMain}.
\end{Remarks*}

%s6 #&#
\section{Discussion}

In this paper we have presented a unified framework for deriving
error bounds and convergence rates for a class of regularized
$M$-estimators. The theory is high-dimensional and nonasymptotic in
nature, meaning that it yields explicit bounds that hold with high
probability for finite sample sizes and reveals the dependence on
dimension and other structural parameters of the model. Two
properties of the $M$-estimator play a central role in our framework.
We isolated the notion of a regularizer being \textit{decomposable} with
respect to a pair of subspaces and showed how it constrains the error
vector---meaning the difference between any solution and the nominal
parameter---to lie within a very specific set. This fact is
significant, because it allows for a fruitful notion of
\textit{restricted strong convexity} to be developed for the loss
function. Since the usual form of strong convexity cannot hold under
high-dimensional scaling, this interaction between the decomposable
regularizer and the loss function is essential.

Our main result (Theorem~\ref{ThmMain}) provides a deterministic bound
on the error for a broad class of regularized $M$-estimators. By
specializing this result to different statistical models, we derived
various explicit convergence rates for different estimators, including
some known results and a range of novel results. We derived
convergence rates for sparse linear models, both under exact and
approximate sparsity assumptions, and these results have been shown to
be minimax optimal \cite{RasWaiYu09}. In the case of sparse group
regularization, we established a novel upper bound of the oracle type,
with a separation between the approximation and estimation error
terms. For matrix estimation, the framework described here has been
used to derive bounds on the Frobenius error that are known to be
minimax-optimal, both for multitask regression and autoregressive
estimation~\cite{NegWai09}, as well as the matrix completion
problem~\cite{NegWai10b}. In recent work \cite{AgaNegWai11}, this
framework has also been applied to obtain minimax-optimal rates for
noisy matrix decomposition, which involves using a combination of the
nuclear norm and elementwise $\ell_1$-norm. Finally, as shown in the
paper \cite{Neg09}, these results may be applied to derive convergence
rates for generalized linear models. Doing so requires leveraging
that restricted strong convexity can also be shown to hold for these
models, as stated in the bound (\ref{EqnSahandGLM}).

There are a variety of interesting open questions associated with our
work. In this paper, for simplicity of exposition, we have specified
the regularization parameter in terms of the dual norm ${\Regplain^*}$ of
the regularizer. In many cases, this choice leads to optimal
convergence rates, including linear regression over $\ell_\qpar$-balls
(Corollary~\ref{CorLassoWeak}) for sufficiently small radii, and
various instances of low-rank matrix regression. In other cases, some
refinements of our convergence rates are possible; for instance, for
the special case of linear sparsity regression (i.e., an exactly
sparse vector, with a constant fraction of nonzero elements), our
rates can be sharpened by a more careful analysis of the noise term,
which allows for a slightly smaller choice of the regularization
parameter. Similarly, there are other nonparametric settings in
which a more delicate choice of the regularization parameter is
required \cite{KolYua10,RasWaiYu10b}. Last, we suspect that there are
many other statistical models, not discussed in this paper, for which
this framework can yield useful results. Some examples include
different types of hierarchical regularizers and/or overlapping group
regularizers \cite{Jac09,JenMai10}, as well as methods using
combinations of decomposable regularizers, such as the fused
Lasso \cite{Tib05}.

\section*{Acknowledgments}

All authors were partially supported by NSF\break Grants DMS-06-05165 and
DMS-09-07632. B. Yu acknowledges additional support from NSF
Grant
SES-0835531 (CDI); M. J. Wainwright and S. N. Negahban acknowledge
additional support from the NSF Grant CDI-0941742 and AFOSR Grant\break
09NL184; and P.~Ravikumar acknowledges additional support from NSF Grant
IIS-101842. We thank a number of people, including Arash Amini, Francis
Bach, Peter Buhlmann, Garvesh Raskutti, Alexandre Tsybakov, Sara van de
Geer and Tong Zhang for helpful discussions.

\begin{supplement}%[id=suppA]
\stitle{Supplementary material for ``A unified framework for high-dimensional
analysis of $\bolds{M}$-estimators with decomposable regularizers''\\}
\slink[doi]{10.1214/12-STS400SUPP} %[doi,text={...}] - jei reikia
%suskaldyti doi
\sdatatype{.pdf}
\sfilename{sts400\_supp.pdf}
\sdescription{Due to space constraints, the
proofs and technical details have been given in the supplementary
document by Negahban et al. \cite{Neg12supp}.}
\end{supplement}

% imsref loaded by lrinkeviciute, 2012-09-26 16:18:54
% imsref loaded by lrinkeviciute, 2012-09-27 08:19:50
% imsref loaded by lrinkeviciute, 2012-09-27 08:43:44

\end{document}